\documentclass[11pt]{amsart}
\usepackage{amssymb,amsfonts,amscd,verbatim}
\usepackage[leqno]{amsmath}
\usepackage{mathrsfs}
\usepackage{xy} 

 \swapnumbers

\numberwithin{equation}{subsection}
\newtheorem{thm}[equation]{Theorem} 
\newtheorem{propose}[equation]{Proposition}
\newtheorem{lemma}[equation]{Lemma}
\newtheorem{cor}[equation]{Corollary}
\theoremstyle{definition}
\newtheorem{defn}[equation]{Definition}
\newtheorem{remark}[equation]{Remark}

\newtheorem{example}[equation]{Example}
\newtheorem{examples}[equation]{Examples}

\renewcommand{\d}{\mbox{\LARGE $\cdot $}}
\newcommand{\Xs}{X_{\d}}            
\renewcommand{\hat}{\widehat}
\newcommand{\Lie}{{\rm Lie}\,}       
\newcommand{\Spec}{\operatorname{Spec}} 
\newcommand{\Spf}{\operatorname{Spf}}   
\newcommand{\Hom}{\operatorname{Hom}}      
\newcommand{\homo}[2]{\mathrm {Hom}(#1,#2)} 
\newcommand{\Ext}{\operatorname{Ext}}      
\newcommand{\eext}[2]{{\mathrm {Ext }}(#1,#2)} 
\newcommand{\uext}[2]{\underline{\mathrm{Ext}}(#1,#2)} 
\newcommand{\Biext}{\operatorname{Biext}}      
\newcommand{\biext}[3]{{\mathrm {Biext}}^1(#1,#2;#3)} 
\newcommand{\mG}{{G}_{\times} } 
\newcommand{\DM}{{\rm DM}}          
\newcommand{\M}{\mathcal{M}_1}   
\newcommand{\Ma}{\mathcal{M}_1^{\rm a}} 
\newcommand{\eMa}{\mathcal{M}_1^{\rm a, eff}} 
\newcommand{\Mlau}{\mathcal{M}_1^{\rm a,fr}} 
\newcommand{\Ml}{\mathcal{M}_{1}^{\ell}} 

\newcommand{\FHSf}{\mathrm{FHS}_{1}^{\mathrm{fr}}} 
\newcommand{\FHS}{\mathrm{FHS}_{1} } 
\newcommand{\MHS}{\mathrm{MHS}_{1}} 
\newcommand{\Tint}{T_{\oint_{}} }

\newcommand{\Div}{\operatorname{Div}}

\newcommand{\ihom}{{\rm\underline{Hom}}}  

\newcommand{\C}{\mathbb{C}}     
\newcommand{\Z}{\mathbb{Z}}     

\newcommand{\G}{\mathbb{G}}     
\newcommand{\HH}{\mathbb{H}}    
\newcommand{\bPic}{{\rm\mathbb{P}ic}}    

\newcommand{\im}{\operatorname{Im}}        
\renewcommand{\ker}{\operatorname{Ker}}  
\newcommand{\coker}{\operatorname{Coker}} 
\newcommand{\Pic}{{\rm Pic}}     
\newcommand{\Alb}{{\rm Alb}}     

\newcommand{\by}[1]{\stackrel{#1}{\rightarrow}}
\newcommand{\longby}[1]{\stackrel{#1}{\longrightarrow}}

\newcommand{\eprooff}{\hfill$\Box$ \\ }

\renewcommand{\tilde}{\widetilde}
\newcommand{\df}{\mbox{\,${:=}$}\,}
\newcommand{\ie}{{\it i.e.\/},\ }
\newcommand{\cf}{{\it cf.\/}\ }
\newcommand{\eg}{{\it e.g.\/},\ }
\newcommand{\id}{{\rm id}}
\newcommand{\et} {{\mbox{\scriptsize{\rm {\'e}t}}}}
\newcommand{\loccit}{{\it loc. cit. }}
\newcommand{\fr}{\mbox{\scriptsize{\rm fr}}}
\newcommand{\tor}{\mbox{\scriptsize{\rm tor}}}
\newcommand{\tf}{\mbox{\scriptsize{\rm tf}}}
\newcommand{\eff}{\mbox{\scriptsize{\rm eff}}}
\newcommand{\sing}{\mbox{\scriptsize{\rm sing}}}
\newcommand{\crys}{\mbox{\scriptsize{\rm crys}}}
\newcommand{\dR}{\mbox{\scriptsize{\rm dR}}}

\newcommand{\mgr}{{\mathbb G}_{m}} 
\newcommand{\agr}{{\mathbb G}_{a}} 
\newcommand{\fagr}{{\hat{\mathbb G}}_{a}} 
\newcommand{\fmot}{\vec{{\mathbb G}}_{a}} 

\renewcommand{\bar}{\overline}
\newcommand{\into}{\hookrightarrow}

\newcommand{\sZ}{\mbox{\scriptsize{$\Z$}}}   

\def\bm#1{\mathpalette\bmstyle{#1}} 
\def\bmstyle#1#2{\mbox{\boldmath$#1#2$}}
\newcommand{\nobd}{\nobreakdash}

\newcommand{\sC}{\mbox{\scriptsize{$\C$}}}   

\newcommand{\liminv}[1]{\mathop{\rm
lim}_{\buildrel\longleftarrow\over{#1}}}
\newcommand{\onto}{\mbox{$\to\!\!\!\!\to$}}

\newcommand{\boxtensor}{\def\boxtimesten{\Box\kern-7.59pt\raise1.2pt
\hbox{$\times$} }}                                  

\newcounter{elno}                   

\newcommand{\cL}{\mathcal{L}}

\newcommand{\cO}{\mathcal{O}}

\newcommand{\cV}{\mathcal{V}}

\begin{document}
\input xy     
\xyoption{all}

\title{Sharp de Rham realization}
\author{Luca Barbieri-Viale}
\author{Alessandra Bertapelle}
\address{Dipartimento di Matematica Pura ed Applicata, Universit\`a  degli Studi di
Padova\\Via G. Belzoni, 7\\Padova -- I-35131\\ Italy}
\email{barbieri@math.unipd.it}
\email{bertapel@math.unipd.it}
\begin{abstract} We introduce the {\it sharp}\, (universal) extension of a 1-motive (with additive factors and torsion) over a field of characteristic zero. We define the sharp de Rham realization $T_{\sharp}$ by  passing to the $\Lie$-algebra. Over the complex numbers we then show a (sharp  de Rham) comparison theorem in the category of formal Hodge structures. For a free 1-motive along with its Cartier dual we get a canonical connection on their sharp extensions yielding a perfect  pairing on sharp realizations. We thus provide one-dimensional sharp de Rham cohomology ${\mathrm H}^1_{\sharp-\dR}$ of algebraic varieties.
\end{abstract}

\maketitle
\section*{Introduction}
Grothendieck's idea of  $\natural$-extensions has been largely  employed and exploited to various extents. In Deligne's  construction (see \cite{D}, 10.1.7,  \cf \cite{MM}) for any Deligne 1-motive $M$ over a field $k$,  one obtains a universal  $\G_a$-extension $M^{\natural}$ of $M$ by the vector group $\Ext (M, \G_a)^{\vee}$. This $M^{\natural}$ is a complex of $k$-group schemes which is no more a Deligne 1-motive. Passing to $\Lie$ algebras Deligne defined in \loccit the de Rham realization $T_{\dR}$ according with Grothendieck's description of one-dimensional de Rham cohomology of abelian varieties in characteristic zero and crystalline cohomology in positive characteristics. Actually, pushing these techniques much  further, we can describe ${\mathrm H}^1_{\dR}$ of any (arbitrarily  singular) algebraic scheme, as well as the (de Jong) ${\mathrm H}^1_{\crys}$ in positive characteristics, by means of universal $\G_a$-extensions of certain Picard 1-motives (and more, \eg see \cite{THS} for the full picture). 
 
\subsection{Results}  We here deal with a {\it sharp}\, $\G_a$-extension, based on the named universal $\G_a$-extension but including $\G_a$-factors, \ie  for Laumon 1-motives, and we show a ``{\it new}\, ${\mathrm H}^1$" out of it. We actually work in the abelian category $\Ma$ of 1-motives with torsion and additive factors, containing the category of Laumon 1-motives $\Mlau$ as the Quillen exact subcategory of free objects  (see \S 1 below for the definition, \cf \cite{LAU}, \S 5 and \cite{BRS}, \S 1).  

\subsubsection*{ } Firstly, we get the sharp $\G_a$-extension $M^{\sharp} = [F\by{u^\sharp}G^{\sharp}]$ of a (effective Laumon) $k$-1-motive $M = [F  \by{u} G]$, over a field $k$ of characteristic zero, roughly as follows. We show that (see \ref{cnsue} and \ref{iffexistence} below) as soon as $\Hom (M,\G_a) =0$ the universal $\G_a$-extension $M^{\natural}$ exists. 
Denote by $M_{\times}$ the quotient of $M$ by the (maximal) $\G_a$\nobd-factor $V(G)$, \cf  \eqref{vbasicext}. The universal $\G_a$-extension $M^{\natural}_{\times}=[F\by{}  \mG^{\natural}] $  of $M_{\times}$ exists, see \ref{pro.univtimes},  and we have that if $M$ admits a universal $\G_a$-extension then $M^{\natural}=M^{\natural}_{\times}$.
Note that, if $M$ is free and $M^*=[ F'\by{} G']$ is the Cartier dual, then the $\Lie$ algebra $\Lie F'^0$ of the connected formal group is dual to the underlying $k$\nobd-vector space $V(G)$. Moreover, if $M$ admits a universal extension then $M^\natural=M_\times^\natural$ is the Cartier dual of 
$[ F'_\et\times \hat{G'}\to G' ]$ where $F'_\et$ denotes the \'etale quotient and $\hat{G'}$  the formal completion at the origin, see \ref{pro.univext}. Thus set $M^{\sharp}$ to be the Cartier dual of  $\vec  M^*\df [F'\times \hat G'\to G']$ and obtain $M^{\sharp}$ as an extension of $M$ by $\Ext (M_{\times}, \G_a)^{\vee}$, see \eqref{diasharp}: this $M^{\sharp}$ is just the pull-back of $M^{\natural}_{\times}$ along $M\to M_{\times}$ (see \S \ref{sharpue} for details). 

\subsubsection*{ } Secondly  {\it via}\,  $M^{\sharp}$ we obtain $T_{\sharp}$ passing to the $\Lie$ algebra, \ie see \S \ref{sharpmot},  define the {\it  sharp de Rham}\, realization $T_{\sharp}(M)\df \Lie (G^{\sharp})$ in such a way that $T_{\sharp}$ can be extended to an {\it exact}\, functor from the abelian category $\Ma$ to filtered $k$-vector spaces. We also show that $T_{\sharp}$ is compatible with Cartier duality, see \S 5. In fact, for $M$ free with Cartier dual $M^{\prime}$ and  ${\mathcal
P}\in \Biext (M, M^{*};\G_m )$ the Poincar\'e biextension, we get a canonical connection $\nabla^{\sharp}$ on the pull-back ${\mathcal P}^\sharp\in \Biext (M^\sharp, M^{*\sharp};\G_m )$ of ${\mathcal P}$ 
 yielding a {\it perfect}\,  pairing $T_\sharp(M)\otimes T_\sharp(M^{*})\to k$, see \ref{thm.perfM}. Over $k =\C$ we then show a (de Rham)  {\it comparison theorem}, see \ref{compthm}, saying that $T_\sharp (M)$ is the underlying $\C$-vector space of  a formal Hodge structure associated to the ${\sharp}$\nobd-extension, \ie we provide a general formula $T_{\oint} (M^{\sharp})\cong T_{\oint} (M)^{\sharp}  $, by making the sharp construction working in the category of formal Hodge structures ${\rm FHS}_1$ and applying a suitable Hodge realization $T_{\oint}$ (see \S 4 below and \cf \cite{FHS}).

\subsubsection*{ }  Thirdly and finally, see \S 6, the resulting ${\mathrm H}^1$'s of an algebraic $k$-scheme $X$ can be visualized  {\it via}\,  a Laumon 1-motive $\Pic^+_{a}(X)$ and, dually, the ${\mathrm H}_1$'s by  its Cartier dual $\Alb^-_{a}(X)$,  see \cite{SAP} for details. For $k=\C$ we get the ${\mathrm H}^1_{\sharp-\sing}(X)\df T_{\oint}(\Pic^+_{a}(X))$.  For $k$ of characteristic zero we now can set ${\mathrm H}^1_{\sharp-\dR}(X)\df T_{\sharp}(\Pic^+_{a}(X))$. Note that here we have set $\Pic^+_{a}(X)\df [0\to \Pic^0 (X)]$ for $X$ proper over $k$. If $X$ is not proper, say smooth for simplicity, $\Pic^+_{a}(X)$ is given by $ [F \to \Pic^0 (\bar X)]$ where $\bar X$ is a suitable compactification, $X = \bar X -Y$, the \'etale part $F_{\et}$ of the formal group $F$ is given by algebraically equivalent to zero divisors on $\bar X$ supported on $Y$ and $F^0$ has $\Lie$ algebra ${\mathrm H}^1_Y(\bar X, \cO_{\bar X})$ modulo the image of ${\mathrm H}^0(X,\cO_X)$: note that the Cartier dual $\Alb^-_{a}(X)$ is the maximal Faltings-W\"ustholz \cite{FW} additive extension of the Serre's Albanese semi-abelian variety of $X$. Over $k=\C$, in \S 6.2 we also link ${\mathrm H}^1_{\sharp-\dR}(X)$ to the enriched Hodge structures of Bloch and Srinivas \cite{BSE}.

\subsection{Perspectives} First of all we expect that a similar construction holds in positive characteristics providing a $T_{\dag}$ and ${\mathrm H}^1_{\sharp-{\crys}}(X) \df T_{\dag}(\Pic^+_{a}(X))$ will be the sharp crystalline cohomology.
For $k =\C$, the general plan is to associate to an algebraic $\C$-scheme $X$ a formal Hodge structure called ``sharp" singular cohomology ${\mathrm H}^*_{\sharp-\sing}(X)$  which contains, in the underlying algebraic structure, a formal group which is an extension of ordinary singular cohomology, \ie ${\mathrm H}^*_{\sharp-\sing}(X)_{\et} = {\mathrm H}^*(X_{\rm an}, \Z)$. There will be ``sharp" versions of de Rham and crystalline as well as comparison theorems between them.  Moreover, the largest 1-motivic part of ${\mathrm H}^{1+i}_{\sharp-\sing}(X)$, ${\mathrm H}^{1+i}_{\sharp-\dR}(X)$, etc. should be exactly that obtained by applying $T_{\oint}$, $T_{\sharp}$, etc. to an algebraically defined (effective Laumon) 1-motive $\Pic^+_{a}(X, i )$ for $i\geq 0$ (generalizing Deligne's conjectures to 1-motives with additive factors, \cf \cite{THS}). 
The main goal of this paper is the first step in drawing such a picture for ${\mathrm H}^1$'s of these forthcoming ``sharp" cohomologies in zero characteristic. 

\subsection{Notations}
$k$ is a field of characteristic $0$ and $\bar k$ is an algebraic
closure. {\bf Aff}$/k$ is the category of affine $k$\nobd-schemes  endowed
with the fppf topology. {\bf Ab}$/k$ is the category of sheaves of
abelian groups on {\bf Aff}$/k$. Given a free $k$\nobd-module $\mathcal
E$ we denote by $\mathcal E^\vee$ the dual $\Hom(\mathcal E,k)$. For a $k$-vector group $V$ we sometimes denote by $V$ also its Lie algebra. 
Given an algebraic $k$\nobd-group $G$ we denote by $\omega_G$ the $k$\nobd-module of invariant differentials on $G$ and by $\bm \omega_G$ the associated sheaf as well as the associated $k$\nobd-vector group. $\hat G$ will denote the formal completion at the origin of $G$ and $\iota\colon \hat G \into G$ the inclusion. Given a formal $k$\nobd-group $F$ we denote by $F^*$ its Cartier dual.

\section{Laumon $1$-motives}
\subsection{Free $1$-motives}
Recall that a  \emph{Laumon $k$\nobd-$1$\nobd-motive} or a \emph{free
$k$\nobd-$1$\nobd-motive} $M=[u\colon F\to G]$ is a two terms complex (in
degree -1, 0) where $F$ is a formal $k$\nobd-group without torsion, $G$ is
a connected algebraic $k$\nobd-group and $u$ is a morphism in {\bf
Ab}$/k$ (\cf \cite{LAU}, 5.1.1). It is known that any formal
$k$\nobd-group $F$ splits canonically as product $F^0\times F_\et$ where
$F^0$ is the identity component of $F$ and is a connected formal
$k$\nobd-group, and $F_\et=F/F^0$ is  \'etale. Moreover, $F_\et$ admits  a
maximal subgroup scheme $F_{\tor}$, \'etale and finite, such
that the quotient $F_\et/F_{\tor}=F_{\fr}$ is constant of
the type $\Z^r$ over $\bar k$. One says that $F$ is without torsion
if $F_{\tor }=0$. The group $G$ is extension of an abelian
variety $A$ by a linear $k$\nobd-group $L$ that is product of its maximal
subtorus $T$ with a vector $k$\nobd-group $V(G)$.

Morphisms of Laumon $k$\nobd-$1$\nobd-motives are morphisms as  complexes.
We will denote by $\Mlau$ the category of Laumon  $k$\nobd-$1$\nobd-motives.

\begin{propose}
The canonical functor $\Mlau \to D^b({\bf Ab}/k)$ is a
full embedding into the derived category of bounded complexes of
sheaves for the fppf topology on {\bf Aff}$/k$.
\end{propose}
\proof
The proof of this
fact for Deligne $1$\nobd-motives in \cite{RA}, 2.3.1, works also for
Laumon $1$\nobd-motives. Indeed  the vanishings
$\Hom(G,F)=0=\Ext_{{\bf Ab}/k}(G,F)$ still hold because of \ref{homtrivfg} \&
\ref{exttrivfg}.
\eprooff

\subsection{Cartier duality}
We recall here the definition of the  Cartier dual of a free  $1$\nobd-motive
$M=[u\colon F\to G]$. See also \cite{LAU}, 5.2.2. Denote by $M_A\df M/L$
the $1$\nobd-motive $[\bar u\colon F\to A]$ induced by $M$ {\it via}\,
the projection $G\onto A$.
The Cartier dual of $M$ is the
$1$\nobd-motive $M^*\df [u'\colon F'\to G']$ where
\begin{itemize}
\item $F'$ is the formal $k$\nobd-group Cartier dual of the affine  algebraic
       $k$\nobd-group $L$.
\item $G'$ is the algebraic $k$\nobd-group that
represents the sheaf on ${\bf Ab}/k$ 
\[\uext{M_A}{\G_m } \colon 
\quad S \leadsto
\Ext_{C^{[-1,0]}({\bf Ab}/k)}(M_A,\G_m)= 
   \Hom_{D^b({\bf Ab}/k)}(M_A,\G_m[1]);\]
\item $u'\colon \ihom(L,\G_m)\to \uext{M_A}{\G_m }$ is  the push-out  morphism for
the  sequence
\begin{equation}\label{seq.LMMA} 0\to L\to M \to M_A\to 0.
\end{equation}
\end{itemize}
We spend some words on the representability of $\uext{M_A}{\G_m }$.
Consider the
sequence of Ext sheaves associated to
\begin{equation}\label{seq.AMAF} 0\to A\to  M_A\to F[1]\to 0.
\end{equation}
It provides
\begin{equation}\label{seq.FDGIAI}
0\to \ihom({F},{\G_m})\to \uext{M_A}{\G_m }\to A' \by{\rho} \uext{F}{\G_m }
\end{equation}
where $A'=\uext{A}{\G_m }$ is the dual abelian variety of $A$.
From Lemma~\ref{exttrivfa}    $\rho=0$, the sheaf  $\uext{M_A}{\G_m }$ is extension of $A'$ by an affine  
algebraic $k$\nobd-group and hence representable by an algebraic 
$k$\nobd-group.

\subsection{Exact sequences}\label{cnsdual}
We will say that a sequence of  two terms complexes of fppf-sheaves,  \eg of free  $k$\nobd-$1$\nobd-motives
\begin{eqnarray}\label{exmot}
0\to M_1\by{f} M_2\by{g}
M_3\to 0\end{eqnarray}
 is \emph{strongly exact} if it is exact as sequence of complexes,
\ie for $k$\nobd-$1$\nobd-motives, if the sequence of algebraic  $k$\nobd-groups $0\to G_1\by{f_G}
G_2\by{g_G} G_3\to 0$ is exact as well as the sequence of formal
$k$\nobd-groups $0\to F_1\to F_2\to F_3\to 0$.\footnote{This is the
definition of \emph{exact} sequences in \cite{LAU}.} One can check
that Cartier duality does not preserve in general strongly exact
sequences. This is pointed out in \cite{DCM} for Deligne
$1$\nobd-motives. For example, consider  a non-trivial $l$-torsion point
$a$ of an abelian variety $A$ for  $l$ a prime number. It
corresponds to a $\mgr$-extension $G'$ of $A'$. Moreover $G'$ is
also extension
\begin{equation}\label{seq.counter}
0\to B'\to G'\by{g} \mgr \to 0
\end{equation}
 where $B'$ is an abelian variety
isogenous to $A'$ and the composition of $g$ with the inclusion
$\mgr \to G'$ is the $l$-multiplication. It is immediate to see that
the dual sequence of (\ref{seq.counter}) is not strongly exact (\cf
\cite{DCM}, \S 1.8).

\begin{propose}\label{pro.strongex} 
Let $0\to M_1\by{f} M_2\by{g}
M_3\to 0$ be a strongly exact sequence of free $1$\nobd-motives. Are  equivalent:
\begin{itemize}
\item[$i)$] the dual sequence is strongly exact;
\item[$ii)$] the complex $~\eta_L\colon 0\to L_1\by{f_L} L_2\by{g_L}  L_3\to 0$ is exact;
\item[$iii)$]  the complex $~\eta_A\colon  0\to A_1\by{f_A} A_2\by{g_A}  A_3\to 0$ is exact.
\end{itemize}
\end{propose}
\proof Expand (\ref{exmot}) writing a diagram having the $0\to
L_i\to G_i\to A_i\to 0$  as vertical sequences. The exactness of
$\eta_L$ , (resp. $\eta_A$) may fail only at $L_2$ (resp. at $A_1$).
Indeed the cokernel of $g_L$ is a linear $k$\nobd-group. It is trivial,
because quotient of the kernel of $g_A$. This implies also the
exactness of $\eta_A$ at $A_2$. Now, it is immediate to check that the
exactness of $\eta_L$  at $L_2$ is equivalent to the exactness of
$\eta_A$ at $A_1$. Hence $ii)\Leftrightarrow   iii)$. Furthermore
$i)\Rightarrow ii)$ because the dual sequence of $\eta_L$ is the
sequence of formal groups. Conversely, $iii)$ implies that  the
induced complex  $0\to M_{A_1}\to M_{A_2}\to M_{A_3}\to 0$  is
strongly exact and hence, passing to duals, we get an exact sequence
of algebraic $k$\nobd-groups $0\to G_1'\to G_2'\to G_3'\to 0$; the
complex of formal $k$\nobd-groups   $F_i'$ is exact because of $ii)$.
\eprooff

\begin{remark}\label{rem.strong} If $M_1$ is a linear $k$\nobd-group or  
if $G_3=0$, the dual sequence is always strongly exact. 
\end{remark}

\subsection{$1$-motives with torsion}
The previous section \ref{cnsdual} motivates the introduction of $1$\nobd-motives with
torsion as done in \cite{BRS} for Deligne $1$\nobd-motives.
\begin{defn} An \emph{effective} (Laumon) $k$\nobd-$1$\nobd-motive $M$
 is a two terms complex (in degree $-1, 0$) $[u\colon F\to G]$ where
$F$ is a formal $k$\nobd-group, $G$ is a connected algebraic  $k$\nobd-group
and $u$ is a morphism in {\bf Ab}$/k$. An \emph{effective} morphism
\[f=(f_F,f_G)\colon [F_1\by{u_1} G_1]\to [F_2\by{u_2} G_2]\] is a
morphism of complexes. The corresponding category is denoted by $\eMa$.

An effective $1$\nobd-motive $M$ is said to be \emph{\'etale} (resp.
\emph{connected}, resp. \emph{special}) if $F^0=0=V(G)$ (resp.
$F=F^0$ and $A=0=T$, resp. $F^0$ maps to $V(G)$ {\it via}\, $u$).  \'Etale
(resp. connected, resp. special) $1$\nobd-motives are regarded as a full
subcategory of $\eMa$.
\end{defn}

Note that since the base field $k$ is assumed of zero characteristic we  always have that $F_{\tor}\times_G V(G) = 0$. In particular, if $u$ is an isomorphism  then $M=0$.

Further, an effective map $f$ is a \emph{quasi-isomorphism}, \ie
$\ker(u_1)\cong\ker(u_2)$ and $\coker(u_1)\cong\coker(u_2)$, if and  only if  $f_G$
is an isogeny and $F_1\cong G_1\times_{G_2} F_2$; see \cite{BRS} for  the more
classical \'etale case. Note that, since char$(k) =0$, we have that  the isogeny $f_G: G_1\to G_2$ is given by pulling back isogenies  between their semi-abelian quotients; in fact, it is an isomorphism  when restricted to the vector groups $V(G_1)\by{\simeq} V (G_2)$.

For an effective $1$\nobd-motive $M$ let  
\[F_{\tor}\cap \ker(u)\df F_{\tor}\times_F \ker(u)\quad {\mathrm{and }}\quad  u(F_{\tor})\df F_{\tor}/ F_{\tor}\cap \ker(u)\subseteq G. \]
We denote
\begin{itemize}
\item $M_{\tor}:=[F_{\tor}\cap \ker(u)\to 0]$, the torsion part of $M$;
\item $M_{\fr} :=[F/F_{\tor}\to G/u(F_{\tor})]$, the free part of $M$;
\item $M_{\tf} :=[F/F_{\tor}\cap \ker(u) \to G] $ the torsion free  part.\\
\end{itemize}
Note that all these operations leave untouched $F^0$ and $V (G)$.
One says that $M$ is \emph{torsion free} if $M_{\tor}=0$; this does
not imply that $F$ is torsion free! Note that $M \to M_{\fr}$ factors as
\[M \by{\psi} M_{\tf} \by{\phi} M_{\fr}\]
where $\phi$ is a quasi-isomorphism and we always have a strongly exact  sequence
\[0\to M_{\tor}\to M\by{\psi} M_{\tf}\to 0
\]
where $\psi$ is a strict epi-morphism. Recall that here, similarly to  \cite{BRS}, \S 1, we say that an  effective morphism $f$ is  \emph{strict} if $f_G$ has connected kernel.
\begin{lemma}[\cf \protect{\cite{BRS}, Prop. 1.3}]
 Let $f=(f_F,f_G):M_1\to M_2$ be an effective  morphism. Then $f$ can  be factored as follows
\begin{equation}\label{strict}
\begin{array}{c}
M_1\longby{f} M_2\\
\scriptstyle \tilde f\displaystyle
\searrow\hspace*{0.5cm} \nearrow\\
\tilde M_2
\end{array}
\end{equation}
where  $\tilde f$ is a strict morphism and
$\tilde M_2 \to M_2$ is a quasi-isomorphism. \end{lemma}
\proof  It follows from \loccit  by pulling back the corresponding
isogeny of the semi-abelian scheme quotients.
\eprooff

Furthermore the class of quasi-isomorphisms admits calculus of right  fractions (\cf \cite{BRS}, Prop. 1.2, and \cite{DCM}, Appendix~C). Define  \emph{morphisms} of 1-motives by inverting quasi-isomorphisms from the  right, \ie  a morphism is represented by
$ f g^{-1} $ where $f$ is effective and $ g$ is a quasi-isomorphism.
\begin{defn} \label{defaddspec} Denote by $\Ma$ (resp.  $\M^s$) the  category of 1-motives \emph{with torsion} obtained localizing the  category of effective Laumon (resp.  special) $1$\nobd-motives  at the  multiplicative class of quasi-isomorphisms.
\end{defn}

\begin{propose} $\Ma$ is an abelian category. $\Mlau\subset \Ma$ is a  Quillen exact sub-category such that $M\leadsto M_{\fr}$ is s  left-adjoint to the embedding.
\end{propose}
\proof Since \eqref{strict} is granted the proof is similar to \cite{BRS}, \S 1, and the more detailed Appendix C in \cite{DCM}.\eprooff

Note that $\Mlau\subset \Ma$ is providing $\Mlau$ of an exact structure in such a way  that the dual of (\ref{seq.counter}) is exact. More generally, Cartier  duality is exact but we won't make use of this fact so that we  omit the proof. (One reduces itself easily to the \'etale case and \cf \cite{DCM}, \S 1.8.)

In the following we denote by $\M$ the category of \'etale  $1$\nobd-motives with torsion regarded as a full abelian sub-category  of $\Ma$ (remark that being \'etale is preserved by quasi-isomorphims).  Deligne 1-motives $\M^{\fr}$ are free \'etale $1$\nobd-motives.
Similarly $\M^s \subset \Ma$ is the full abelian sub-category of  1-motives with torsion that are special.

\subsection{Linearized $1$-motives}
For the sake of exposition we introduce a category $\Ml$ which is equivalent to $\Ma$ but 
where connected formal groups are replaced by $k$-vector spaces, since char$(k) =0$.

\begin{defn} Let $\M^{\ell , \eff}$ be the category whose objects are  pairs
$(u_\et\colon  F_\et\to G, u_a\colon F_a\to \Lie(G))$
with $G$ a connected algebraic $k$\nobd-group, $F_\et$ an \'etale formal $k$-group, $u_\et$  a morphism in  {\bf Ab}$/k$
and $u_a$ a homomorphism of finite dimensional $k$\nobd-vector spaces.
Effective morphisms are  triples
\[(f_\et\colon F_\et\to F_\et' , f\colon G\to G', f_a\colon  F_a\to  F_a'  )\]
with $f_\et,f$ morphisms in ${\bf Ab}/k$ and  $f_a$ a homomorphism of vector spaces such that
 the obvious squares commute. Let $\Ml$ be the category obtained by localizing from the right $\M^{\ell  , \eff}$ at the multiplicative class of quasi-isomorphisms on the first component (and isomorphisms on the second).
\end{defn}

\begin{propose} \label{linear}
Let $[u\colon F\to G] $ be a $1$\nobd-motive.
The functor
\[\Ma\to \M^\ell , \quad [u\colon F\to G] \mapsto
\big( u_\et\colon  F_\et\to G,u_a\colon \Lie( F^0) \longby{ \Lie(u^0)
}\Lie(G)\big)\]
is an equivalence of categories.
\end{propose}
\proof
Given a pair   $( u_\et\colon  F_\et\to G,
u_a\colon F_a\to \Lie(G))$ as above we get a connected  formal
$k$\nobd-group $F^0$ as the formal completion at the origin of the  vector
group $\Spec(k[F_a^\vee])$. Moreover, as $\hat G$ is isomorphic to
the formal completion at the origin of  $\Spec(k[\Lie(G)^\vee])$
(\cf \ref{subformcompl})
 the homomorphism $u_a$ provides a morphism of formal
$k$\nobd-groups $F^0\to \hat G$ and  hence  a morphism $F^0\to G$ in   ${\bf Ab}/k$.
\eprooff

The category $\M^\ell$ is somewhat meaningful in order to construct objects in $\Ma$
from geometric invariants associated to algebraic schemes (see \cite{FW} and \cite{SAP}).

\section{Universal extensions of  $1$-motives}
Let $M = [u\colon F \to G]$ be an effective $k$\nobd-$1$\nobd-motive.

\subsection{Some notations}\label{uenot} For $M = [F \by{u} G]$ over $k$ let
$V (G)\subseteq G$ be the maximal vector subgroup of $G$ so that $G$
can be represented as  follows
\begin{equation}\label{gmext} 0\to V (G)\to G \to \mG \to 0
\end{equation} where $\mG$ is the semi-abelian quotient and
$V(G)\cong \G_a^n$ for some $n$.
Denote $u_{\times}$ the composition of $u$ and the projection $G\onto  \mG$. Set
\begin{itemize}
\item $M_{\times}\df [u_{\times}\colon F\to \mG ]$
\end{itemize}
in such a way that we have a short exact sequence of complexes
\begin{equation}\label{vbasicext} 
0\to V (G)[0]\to M \to M_{\times}\to 0.
\end{equation}
Moreover, recalling that $F=F^0\times_{k}F_{\et}$ canonically, we
denote by $u_{\et}$ the composition of  $F_{\et}\into F$ and
$u_{\times}$. Set
\begin{itemize}
\item $M_{\et}\df [u_{\et}\colon F_{\et}\to \mG ].$
\end{itemize}
It is an \'etale $1$\nobd-motive; if $F_{\et}$ is free, $M_{\et}$ is a  Deligne $1$\nobd-motive. We further get a functor 
$M\leadsto M_{\et}:  \Ma\to \M$, left inverse to the inclusion of \'etale 1-motives. We  always have a strongly exact sequence
\begin{equation}\label{fbasicext} 
0\to M_{\et}\to M_{\times}\to  F^0[1]\to 0.
\end{equation}
Given a connected algebraic $k$-group $G$ denote $\vec G\df [\hat {G}\by{\iota} G]$ the induced $1$\nobd-motive. Set moreover 
\begin{itemize}
\item $\vec M\df  [\vec u \colon F\times \widehat{G}\to G ]$
\end{itemize}
obtained as push-out of $M$ with respect to $G\to \vec G$. We also have the restriction of $\vec M$ to $\widehat{V(G)}$, \ie  $[F\times \widehat{ V(G)}\to G ]$ which is also an extension of $M_{\times}$ by $\overrightarrow{V(G)}$.
If $M$ is special, we can further set $M^0 \df [F^0\to V (G)]$, and  then get
\begin{equation}\label{fspecext} 
0\to M^0\to M\to M_{\et}\to 0
\end{equation} so that $M\leadsto M_{\et}: \M^s\to \M$ is left adjoint  to the inclusion $\M \into \M^s$ (\cf \cite[\S 2]{FHS}).

\begin{defn} Let $M$ be an effective $k$\nobd-$1$\nobd-motive such that
$\Hom(M,W)=0$ for any $k$\nobd-vector group $W$. We say that $M$ admits  a \emph{universal
$\G_a$-extension} if it exists a $k$\nobd-vector group ${\mathbb V}(M)$
and an extension $M^\natural$ of $M$ by ${\mathbb V}(M)$ such that
the push-out homomorphism
 \begin{equation}\label{defunivext}
 \Hom({\mathbb V}(M),W   )\longrightarrow \Ext(M,W)
\end{equation}
is an isomorphism for any $k$\nobd-vector group $W$.
 It is immediate to
check that ${\mathbb V}(M)$ has to be then the vector group
associated to $\Ext(M,\G_a)^\vee$.
\end{defn}

 Observe that thanks to \ref{lem.ext} \& \ref{pro.ext}
the notation $\Ext(M,W)$ is not ambiguous.

\begin{remark}\label{cnsue} Following \cite{MM}, I, 1.7,
one can see that $M$ admits a universal extension
$M^\natural$ if and only if
\begin{itemize}
\item[a)] $\homo{M}{\G_a}=0$,
\item[b)] $\eext{M}{\G_a}$ is a $k$\nobd-vector space of finite  dimension.
\end{itemize}
If $M^\natural$ exists, then  $\Hom(M^\natural,\G_a)=0=\Ext(M^\natural,\G_a)$
and $M^{\natural\natural}=M^\natural$.
\end{remark}
\begin{examples}\label{exue}We have the following paradigmatic cases:
\begin{itemize}
\item  $\G_a$ does not admit universal extension. \item
$T^\natural=T$ for any $k$\nobd-torus $T$. \item For any abelian
variety $A$ the universal extension $A^\natural$ exists (\cf
\cite{MM}). As observed in \cite{LAU}, 5.2.5, $A^\natural$ is the
Cartier dual of the $1$\nobd-motive $\vec A^\prime = [\hat A^\prime\to A^\prime]$.
\item Any Deligne $1$\nobd-motive $M_{\et} = [u_{\et}\colon
F_{\et}\to G_\times]$ ($F_{\et}$ free) admits a universal
$\G_a$-extension $M^{\natural}_{\et} = [u_{\et}^{\natural}\colon
F_{\et}\to G^{\natural}_{\et} ]$ (see \cite{D}).\end{itemize}
\end{examples}

\subsection{Existence results}
We start showing that we can reduce to work with free $1$\nobd-motives.

\begin{propose}\label{pro.uefr}
An effective $1$\nobd-motive $M$ admits universal extension if and only
if $M_{\fr}$ does.
\end{propose}
\proof Set $ K:=F_{\tor}\cap \ker(u)$ and consider the sequence
\begin{eqnarray}\label{seqMMtf}
0\to K[1]=M_{\tor}\to M\to M_{\tf}\to 0.
\end{eqnarray}
As $\Hom(M_{\tor},\G_a)=0=\Ext(M_{\tor},\G_a)$, conditions in  \ref{cnsue}
holds for $M$ if and only if they hold for $M_{\tf}$. Moreover, if
$M_{\tf}=[F/K\by{v} G^\natural]$ is the universal extension of
$M_{\tf}$,  the universal extension of $M$ is $[F\to G^\natural]$
obtained {\it via}\, composition of $v$ with the canonical  $F\to F/K$.
As $M_{\tf}$ and
$M_{\fr}$ are quasi-isomorphic, conditions in  \ref{cnsue} hold for
both or none of them. Moreover if the universal extensions
$M_{\tf}^\natural$ and $M_{\fr}^\natural$ exist, they are
quasi-isomorphic and $M_{\tf}^\natural$ is obtained {\it via}\,  pull-back of
$M_{\fr}^\natural$ along $M_{\tf} \to M_{\fr}$.
\eprooff

 We will see that for effective $1$\nobd-motives, condition $a)$
in \ref{cnsue} implies condition $b)$.
 We start with the case $M=M_\times$.

\begin{propose}\label{pro.univtimes}
Let $M$ be an effective $1$\nobd-motive. The universal $\G_a$-extension
\[M^{\natural}_{\times}=[u_{\times}^{\natural}\colon F\to  \mG^{\natural}] \]
 of $M_{\times}$ exists.
\end{propose}
\proof By \ref{pro.uefr} we may suppose $F$ torsion free. As
$G_\times$ is semi-abelian, $M_\times$ satisfies condition $a)$ in
\ref{cnsue}; moreover, from (\ref{fbasicext}), \ref{lem.ext} \& 
\ref{exttrivfa} we
get
\[0\to \Hom(F^0,\G_a) \to \Ext(M_\times,\G_a) \to   \Ext(M_\et,\G_a)  \to \Ext(F^0,\G_a)=0.
\]
Now, $M_{\et}$ is a Deligne $1$\nobd-motive and hence  $M_\et$
satisfies condition $b)$ (\cf \ref{exue}). As $\Hom(F^0,\G_a)$ is a
free $k$\nobd-module (\cf \ref{lem.homform}) also $M_\times$ satisfies
condition $b)$ and we are done. \eprooff

As for Deligne $1$\nobd-motives, we have the following description of
$\Ext(M_\times,\G_{a})^\vee$ in terms of invariant differentials:

\begin{lemma} \label{Gtimesnat}
Let $M^* =[u' : F'\to G']$ be the Cartier dual of $M$ free. Then
\begin{eqnarray}\label{eq.lievect}
    \Ext(M_\times,\G_{a})^\vee= \Lie(G')^\vee=\omega_{G'}.
\end{eqnarray}
Moreover, $G_\times^\natural$ is the push-out of $A^\natural$ with
respect to the canonical homomorphism
\[\Ext(A',\G_a)^\vee={\bm \omega}_{A'}\longrightarrow {\bm \omega}_{G'}=
\Ext(M_\times,\G_a)^\vee=\Ext(M_A,\G_a)^\vee.\]
\end{lemma}
\proof The arguments in \cite{BER}, 2.6, work also for the effective
$1$\nobd-motive $M_\times$. The second assertion can be checked as
for Deligne $1$\nobd-motives (\cf \cite{BVS}). \eprooff

As $M$ is extension of $M_\times$ by the vector group $V(G)$, we may
view $M$ as the push-out of the universal extension
$M_\times^\natural$ of $M_\times$ with respect to a unique $v_M$:
\begin{eqnarray}\label{diavM}
\begin{CD}
0@>>> \Ext (M_{\times},\G_a)^{\vee}@>{}>> M^{\natural}_{\times}@>{}>>
M_\times @>>> 0\\
&&@V{v_M}VV @V{v_M^\natural}VV@V{||}VV\\
0@>>>V(G)@>>> M @>>> M_\times @>>> 0
\end{CD}\end{eqnarray}

\begin{propose}\label{iffexistence}
Let  $M=[u\colon F\to G]$  be an effective $1$\nobd-motive.   Are
equivalent:
\begin{itemize}
\item[$i)$] $M$ admits a universal $\G_a$-extension,
\item[$ii)$] $\Hom(M,\G_a)=0$,
\item[$iii)$]$v_M$ is surjective.
\end{itemize}
Moreover, if $M$ admits a universal extension then  $M^\natural=M^\natural_\times$.
\end{propose}
\proof We start showing that  $\Hom(M,\G_a)=0$ if and only if  $v_M$ is  surjective.
From (\ref{vbasicext}) we get a sequence
\[0\to \Hom(M,\G_a)\to  \Hom(V(G),\G_a)\by{\partial}  \Ext(M_\times,\G_a)\to \Ext(M,\G_a)\to 0.  \]
Now, $\Hom(M,\G_a)=0$ if and only if $\partial$ is injective.
However, $\partial$ is the push-out homomorphism along $v_{M}$ followed
by the isomorphism
$\Hom(\Ext(M_\times,\G_a)^\vee,\G_a)=\Ext(M_\times,\G_a)$ (\cf
(\ref{defunivext})). Hence $\partial$ is injective if and only if
the push-out along $v_{M}$ is injective and this last is equivalent to the
surjectivity of $v_M$. Furthermore,  $\partial$ coincides with the
$\G_a$-dual map $v_M^\vee$ and hence we get the exact sequence
\begin{eqnarray}\label{seqkervM}
0\to \Ext(M,\G_a)^\vee\to \Ext(M_\times,\G_a)^\vee \by{v_M} V(G)\to 0.
\end{eqnarray}
Moreover $i)\Rightarrow ii)$ by definition.

Suppose now that $v_M$ is surjective (an hence $\Hom(M,\G_a)=0$).
From (\ref{seqkervM}) and \ref{pro.univtimes} we get that
$\Ext(M,\G_a)^\vee=\ker(v_M)$ is a finite-dimensional $k$-vector space.
Hence $M$ clearly satisfies condition
$b)$ and it admits a universal extension.

For the last assertion, note that $\ker(v_M)=\ker(v_M^\natural)$ (\cf
(\ref{diavM})  and (\ref{seqkervM})); it is immediate to check that
\[0\to  \Ext(M,\G_a)^\vee \to M_\times^\natural\by{v_M^\natural}  M\to  0\]
 satisfies the universal property.
\eprooff

\begin{examples} a) $\fmot:=[\iota\colon \fagr\to \G_a]$, with $\iota$ the
inclusion, is the universal extension of $\fagr[1]$. Note that the Cartier dual $\fagr[1]^*=\G_a$ does not admit a universal extension. More generally,
let $F$ be a connected formal $k$\nobd-group. The universal extension of
$F[1]$ is the $1$\nobd-motive $\overrightarrow{\Lie(F)}=[F\by{\iota} \Lie(F)]$. To show this fact, one uses
\[\Ext(F[1],\G_a)^\vee=\Hom(F,\G_a)^\vee  =\Hom(\Lie(F),\G_a)^\vee=\Lie(F).\]
b) Let $F$ be a torsion formal $k$\nobd-group.
Then $F[1]=F[1]^\natural$.
\end{examples}

By making use of Laumon $1$\nobd-motives one can give an interpretation  of universal
extensions in terms of dual $1$\nobd-motives.

\begin{propose}\label{pro.univext}
Let $M=[u\colon F\to G]$ be a free $1$\nobd-motive and $M^*=[u'\colon
F'\to G']$ its Cartier dual. If $M$ admits a universal
extension $M^\natural=[F\to G^\natural]$ then
$M^\natural=M_\times^\natural$ is the Cartier dual of $[(u',\iota)\colon F'_\et\times \hat G'\to G' ]$ where \[
(u',\iota)(x,y)\df  u'(x)+\iota (y).\]
\end{propose}
\proof We have already seen that $M^\natural=M_\times^\natural$. Hence  we
may reduce to the case $M=M_\times$.
We start with the case $M=M_A=[F\to A]$. Let
\[0\to {\bm \omega}_{G'}\to E\to M_A\to 0\]
be the Cartier dual of $\vec G' = [\hat G'\by{\iota} G']$ (\cf \ref{lem.conno}).
Suppose given an extension $N$ of $M_A$ by a $k$\nobd-vector group $W$.
The Cartier dual of $N$ is an extension of $W^*=\ihom(W,\G_m)[1]$
by $G'$ (\cf \ref{rem.strong}), hence it corresponds to a morphism $h\colon \ihom(W,\G_m)\to  G'$.
As $\ihom(W,\G_m)$ is a connected formal $k$\nobd-group, $h$ factors
through a unique morphism $\bar h\colon \ihom(W,\G_m)\to \hat G'$ and
hence $N$ is the push-out of $E$ {\it via}\, the dual morphism  
$ h^*\colon {\bm \omega}_{G'}\to W$.
In particular, $E$ is the
universal $\G_a$-extension of $M_A$.

In the case $T$ is not trivial, $M_\times$ is extension of $M_A$ by $T$.
Its universal extension is the pull-back of $M_A^\natural$ along
$M_\times\to M_A$. Hence the Cartier dual of $M_\times^\natural$
is the push-out of $\vec G'$ along
$G'\to [u'\colon F'_{\et}\to G']=(M_\times)^*$ and this last is the $1$\nobd-motive
$[(u',\iota)\colon F'_\et\times \hat G'\to G' ]$.
\eprooff

\subsection{Exact sequences}\label{exactseque}
It follows from \ref{iffexistence} that given a strongly exact
sequence of $1$\nobd-motives
\begin{equation}\label{eq.exactue}
0\to M_1\to M\to M_2 \to 0\end{equation}
if $M$ admits universal extension, then the same does $M_2$ while
 $M_1$ may not admit universal extension.
For example, consider the sequence
\[0\to \G_a\to \fmot\to \hat \G_a[1]\to 0.\]
Moreover, we have:
\begin{lemma}
If $M_1,M_2$ admit universal extensions also $M$ admits
universal extension and
\[0\to M_1^\natural \to M^\natural\to
M_2^\natural \to 0 \] is strongly exact.
\end{lemma}
\proof
This fact follows immediately
from \ref{pro.univext} if the dual of (\ref{eq.exactue}) is still
strongly exact.
In the general case one has to check that given
 an isogeny of abelian varieties $\varphi\colon A\to B$, a free formal
$k$\nobd-group $F$ and a morphism $u\colon F\to A$, the universal  extension of
$[u\colon F\to A]$ is the pull-back {\it via}\, $\varphi$ of the  universal extension
of $[\varphi \circ u\colon F\to B]$.
\eprooff

\begin{remark}\label{rem.fbasicsh}
In particular, the sequence (\ref{fbasicext}) provides an exact sequence
\[0\to M_\et^\natural\to M_\times^\natural\to
 \overrightarrow{\Lie(F^0)}\to 0.\]
\end{remark}

 \section{Sharp de Rham realization of $1$-motives}\label{sharpue}
\subsection{Sharp (universal) extension}
Proposition \ref{pro.univext} shows that the universal extension, when  it
exists, forgets the contribution of the $k$-vector group $V(G)$, \ie of  the connected  
formal group $F^{'0}$ of the dual.
We introduce then a more general object.

\begin{defn} Let $M=[F\by{u}G]$ be  a free $1$\nobd-motive and
$M^*=[F'\by{u'} G']$ its Cartier
dual. The \emph{sharp $\G_a$-extension} $M^\sharp\df [u^\sharp\colon F\to  G^\sharp]$ of  $M$ is the Cartier dual of the $1$\nobd-motive $\vec  M^*= [(u',\iota)\colon F'\times \hat{G'}\to G']$.
\end{defn}

\begin{lemma}\label{lem.sharp}
The free $1$\nobd-motive $M^\sharp$ fits in the following pull-back
diagram \begin{eqnarray}\label{diasharp}
\xymatrix{
&& V(G)\ar@2{-}[r]\ar[d]& V(G)\ar[d]&\\
 0\ar[r]  &  (\hat G')^*= \Ext
(M_{\times},\G_a)^{\vee}\ar[r] \ar@2{-}[d]& M^{\sharp}\ar[r]\ar[d]  &
M \ar[r]\ar[d]&  0\\
0 \ar[r]& \Ext (M_{\times},\G_a)^{\vee}\ar[r] &
M_{\times}^{\natural}\ar[r]\ar@{.>}[ur]|{v_M^\sharp}&
M_{\times}\ar[r]& 0
}\end{eqnarray}
Moreover, the  homomorphism $v_M^\natural\colon
M_\times^\natural \to M$ of (\ref{diavM}) provides a splitting of
 the  vertical sequence in the middle.
\end{lemma}
\proof
As  $\vec M^*$ fits in the following (strongly exact) sequence
\[0\to M^*\to \vec M^*\to \hat G'[1]\to 0,\]
passing to duals, one gets the horizontal sequence in the middle
of (\ref{diasharp}).
The map $M^\sharp\to M_\times^\natural$ is the dual of
the canonical morphism  $[F'_{\et}\times \hat G'\to G']\to \vec M^*$.
The last assertion is immediate.
\eprooff

\begin{lemma}\label{lem.Gsh}
The algebraic $k$\nobd-group $G^\sharp$ fits in the following diagram
\begin{eqnarray*}
\xymatrix{
 0\ar[r]&{\Ext([F^{\et} \to A],\G_a)^\vee =\bm  \omega}_{A'}\ar[r]\ar@{^{(}->}[d]^i &A^\natural\times_A G \ar[r]  \ar@{^{(}->}[d] &
G\ar[r]\ar@2{-}[d] &0\\
0\ar[r]&\Ext(M_\times,\G_a)^\vee ={\bm \omega}_{G'}  \ar[r]\ar@{->>}[d]^{\bar\tau}& G^\sharp
\ar[r]^\rho\ar@{->>}[d]^\tau &
G\ar[r] &0\\
& {\bm \omega}_{L'} \ar@2{-}[r] & {\bm \omega}_{L'}  & &.
}\end{eqnarray*}
that generalizes the one in \cite{BVS} for Deligne $1$\nobd-motives.
\end{lemma}
\proof
By construction $G_\times^\natural$ is the
push-out of $A^\natural$ with respect to ${\bm
\omega}_{A'}\longrightarrow {\bm \omega}_{G'}$ (see
(\ref{Gtimesnat})) and $G^\sharp$ is the pull-back {\it via}\, $G\to A$  of
$G_\times^\natural$. The previous diagram says that we can take
first the pull-back and then the push-out.
\eprooff

Lemma \ref{lem.sharp} provides an alternative definition of $M^\sharp$ for
free $1$-motives that can be extended to effective $1$-motives.

\begin{defn}
Let $M$ be an effective $1$-motive. Denote by $M^\sharp\df [F\to G^\sharp]$  
the pull-back  of $M_\times^\natural$ along $M\to M_\times$ and call it \emph{sharp} $\G_a$\nobd-\emph{extension} of $M$. In particular \eqref{diasharp} holds.
\end{defn}

By definition, $(M_\times)^\sharp=M_\times^\natural$. 
However this equality does not  hold for a general $M$ that amdits
universal extension. For example, $\fmot^\natural=\fmot$ 
while $\fmot^\sharp=[\hat \G_a\to  \G_a^{2} ]$
with the diagonal embedding as morphism.

\begin{lemma} \label{lem.exsharp}
The functor $(\ \ )^\sharp\colon \Ma\to \Ma$ is exact. 
\end{lemma}
\proof
Suppose given a quasi-isomorphism $f\colon M_1\to M_2$ of effective
$1$\nobd-motives. It induces a quasi-isomorphism $f_\times \colon
M_{1\times}\to M_{2\times}$ and then a quasi-isomorphism
$M_{1\times}^\natural\to M_{2\times}^\natural$.  In particular,
$f^\sharp\colon M_1^\sharp\to M_2^\sharp$ is a quasi-isomorphism.
Any short exact sequence in $\Ma$ is isomorphic to a strongly exact
sequence of effective $1$\nobd-motives. Hence we may restrict to work with 
strongly exact sequences 
of effective $1$-motives as in \ref{eq.exactue}. Arguments used  in 
the proof of \ref{pro.strongex} say that the complex of linear 
$k$\nobd-subgroups is  exact (not necessarily strongly) and it is fixed by the $(\ \ )^\sharp$ functor.
Moreover
\begin{eqnarray}\label{seq.ma}
0\to M_{1,A_1}\to M_A \to M_{2,A_2}\to 0\end{eqnarray}
is exact. We are then reduced  to see that $(\ \ )^\sharp$ applied to 
\eqref{seq.ma} preserves exactness. Thanks to the horizontal sequence
in the middle of \eqref{diasharp}, it is sufficient to check that
\[0\to \Ext(M_{1,A_1},\G_a)^{\vee}\to \Ext
(M_{A},\G_a)^{\vee} \to \Ext
(M_{2,A_2},\G_a)^{\vee}\to 0\]
is exact. Let $B$ be the kernel of 
$A\to A_2$; it is an abelian variety isogenous to $A_1$. 
From \ref{pro.ext} and the isomorphism 
$\Ext([F_1\to B],\G_a)=\Ext(M_{1,A_1},\G_a)$, we get the result.
\eprooff

\subsection{Sharp de Rham}\label{sharpmot} We now can set the following:
\begin{defn}
Let $M$ be an effective $1$\nobd-motive. Its \emph{sharp de Rham realization}  is
$$T_{\sharp}(M) \df \Lie (G^{\sharp}).$$
\end{defn}
Observe that $\Lie (G^{\sharp})$ contains $V (G)$ and  $\Ext
(M_{\times},\G_a)^{\vee}$ in such a way that $\Ext (M_{\et},\G_a)^{\vee}
\subseteq \Ext (M_{\times},\G_a)^{\vee}$ with quotient $\Hom
(F^0,\G_a)^{\vee}$ and we clearly have that $\Ext
(M_{\times},\G_a)^{\vee}\cap V(G) =0$.
The diagram of Lie algebras of (\ref{diasharp}) yields
\[\begin{CD}
&& && 0 &&
0&&\\
&& && @V{}VV@V{}VV\\
&& && V (G) @>{=}>>
V (G) &&\\
&& && @V{}VV@V{}VV\\
0@>>> \Ext (M_{\times},\G_a)^{\vee}@>{}>> \Lie (G^{\sharp})@>{}>>
\Lie (G) @>>> 0\\
&&@V{||}VV @V{}VV@V{}VV\\
0@>>>\Ext (M_{\times},\G_a)^{\vee}@>>>\Lie (G_{\times}^{\natural})@>>>  \Lie (G_{\times})@>>> 0\\
&& && @V{}VV@V{}VV\\
&& && 0 &&
0&&\\
\end{CD}\]
and provides for a free $1$-motive $M$
\[0\to T_{\dR}(M_{\et}) \to T_{\sharp}(M)/V(G)\to \Hom (F^0,  \G_a)^{\vee}\to 0
\]
that is
\[0\to T_{\sharp}(M_{\et}) \to T_{\sharp}(M_{\times})\to T_{\sharp}  (F^0[1])\to 0 .\]
Sharp de Rham realization is compatible with \eqref{vbasicext} and
\eqref{fbasicext}.

\begin{propose} The functor $T_{\sharp}$ behaves well passing to localization on
quasi-isomorphisms and it provides an exact functor from $\Ma$ to the category  of
(filtered) $k$\nobd-vector spaces
\[T_\sharp\colon \Ma \longrightarrow {\cV_k}.\]
\end{propose}
\proof
We have already seen in \ref{lem.exsharp} that $(\ \ )^\sharp$ 
is an exact functor. 
Moreover, any quasi-isomorphism $M_1\to M_2$ induces an isomorphism 
 $\Lie(G_1)\to \Lie(G_2)$.
The conclusion follows recalling that any exact sequence is represented by 
an effective exact sequence.
\eprooff

\begin{remark}\label{rem.sharp}
It follows from the proof of \ref{pro.uefr} that $M^\sharp$ is the  pull-back
of $M^\sharp_{\fr}$ along the canonical morphism $M\to M_{\fr}$ and 
$T_{\sharp}(M) \cong T_{\sharp}(M_{\fr})$.
\end{remark}

\section{Hodge theory}
In this section $k=\C$. Also assume that the mixed Hodge structures are  graded polarizable and denote by $\mathrm{MHS}_1$ the category of those  structures with possibly non-zero Hodge numbers in the set $\{(0,0),  (-1,0),(0,-1),(-1,-1)\}$, \ie of level $\leq 1$. The key point in what  follows (\cf \ref{rem.formlie} and \ref{linear}) is that working with  a connected formal $\C$\nobd-group $F^0$ we  can think of $F^0$ also as  $\Lie(F^0)$, the associated $\C$\nobd-vector group  $\Spec(\C[\Lie(F^0)^\vee])$ or just the underlying $\C$-vector space.

\subsection{Formal Hodge structures}
This section is based on \cite{FHS}.
A formal Hodge structure (of level $\leq 1$) is: $(i)$
 a formal $\C$\nobd-group $H=H^0\times H_\Z$ such that $H_\Z$ admits a mixed Hodge structure $H_{\et}=(H_\Z,W_*,F^*_{Hodge})\in\mathrm{MHS}_1$, $(ii)$ a $\C$\nobd-vector  space $V$ with a
sub-space $V^0\subseteq V$, $(iii)$ a ``group homomorphism'' $v\colon
H\to V$, (\ie a homomorphism of $\C$-vector spaces $v^0\colon \Lie(H^0)\to V$
and a homomorphism of abelian groups $v_\Z\colon H_\Z\to V$) $(iv)$ a
$\C$-isomorphism $\sigma\colon H_\C/F^0_{Hodge}\cong V/V^0$.
Moreover if  $c\colon H_\Z\to H_\C/F^0_{Hodge}$ is the canonical map and 
$pr\colon V\to V/V^0$ is the projection, we assume that the following square
\begin{eqnarray}
\label{dia.fhs}\begin{CD}
H_\Z @>{v_\Z}>> V  \\
  @V{c}VV @V{pr}VV \\
H_\C/F^0_{Hodge}@>{\sigma}>> V/V^0
\end{CD} \end{eqnarray}
 commutes.

 We denote by $V^1$ the sub-space of $V$ that is the pull-back of  $\sigma(W_{-2}H_\C)$ {\it via}\,
 $V\to V/V^0$.   It holds $V^0\subseteq V^1\subseteq V$.
 We denote by $v_\C$ the $\C$-linear map $H_\C\to V$ induced by $v_\Z$.

Denote by $(H,V)$ for short a formal Hodge structure. A morphism between
$(H,V)$ and $(H',V')$ is a morphism of formal groups $h\colon H\to
H'$ and a $\C$-homomorphism $ g\colon V\to V'$ that respects the
above structures and conditions. Denote by $\mathrm{FHS}_{1}$ the  category of
formal Hodge structures and by  $\mathrm{FHS}_{1}^{{\fr}}$ the full
subcategory given by $(H,V)$ with $H$ free. A formal Hodge structure
$(H,V)\in \mathrm{FHS}_{1}$ is said to be \emph{special} (resp.  \emph{connected}) if $v(H^0)$ lies in $V^0$
(resp. $H_\Z=0$). Denote by $\mathrm{FHS}_{1}^{s}$ the full
subcategory of special structures and by   $\mathrm{FHS}_{1}^{0}$
the full subcategory of connected structures. Set
\begin{eqnarray}\label{overtofhs}
\overrightarrow{(H,V)}\df (H\times \hat V, V)\text{\hspace*{0.5cm}and\hspace*{0.5cm}}
\overrightarrow{(H,V)}_0\df (H\times \hat V_0, V)
\end{eqnarray}
where the filtration on $V$ remains the same and the morphism 
$\vec v\colon H\times \hat V\to V$ is induced by $v\colon H\to V$ and $ V  \by{id} V$ (or the inclusion $V_0\subset V$). 

\subsection{Enriched  Hodge structures}\label{sec.ehs} Recall that an enriched
Hodge structure (of level $\leq 1$) is a pair $(H_\et, U\by{u}
V)$ where $H_\et= (H_{\sZ},W_*,F^*_{Hodge})\in \mathrm{MHS}_1$, $u$ is a
$\C$-linear map of $\C$-vector spaces and there exists a commutative
diagram
\begin{eqnarray}\label{dia.ehs1}
\xymatrix @C 1.5cm{
H_\C\ar[r]^{\rho} & U\ar[r]^{\pi}\ar[d]^u & H_{\C}\ar[d]^c \\
& V \ar[r]^{\pi_0}&  H_{\sC}/F^0_{Hodge}
}
\end{eqnarray}
where the composition of the upper arrows is the identity (\cf
\cite{BSE}). The category of such enriched Hodge structures is
denoted by  $\mathrm{EHS}_{1}$. It is clear that we have a functor
\begin{eqnarray}\label{funEF}
\mathrm{EHS}_{1} \longrightarrow  \mathrm{FHS}_{1},
\quad (H_\et, U\by{u} V) \mapsto (H_\Z\times\hat{\ker(\pi)}, V)
\end{eqnarray}
where $\iota\colon \hat{\ker(\pi)}\into \ker(\pi)$ is the completion at the origin  and $v\colon H_\Z\times \ker(\pi)\to  V$ is obtained {\it via}\,  $u$  by composition with the maps $H_\Z\to H_\C\to U$ and $\ker(\pi)\to U$.
The sub-space $V^0$ of $V$ is defined as the kernel of $\pi_0$; hence we
get an isomorphism $\sigma\colon H_{\sC}/F^0_{Hodge}\to V/V^0 $ and the
diagram (\ref{dia.fhs}) commutes by construction.  Furthermore, by
construction $v(\ker(\pi))$ lies in  $V^0$; hence $(H_\Z\times
\hat{\ker(\pi)}, V)$ is special. We actually obtain:
\begin{propose}\label{enrich}
 $\mathrm{EHS}_{1}$ is equivalent to the full subcategory
$ \mathrm{FHS}_{1}^s$ of $ \mathrm{FHS}_{1}$.
\end{propose}
\proof Let $(H,V)$ in $\mathrm{FHS}_{1}$. Note that giving the map 
$v\colon H\to V$ is equivalent to give a map $u\colon \Lie(H^0)\oplus H_\C\to V$ of $\C$\nobd-vector spaces.  If  $(H, V)$ is special then $\Lie (H^0)$ is mapped to $V^0$ and $(H_\et,  H_{\sC}\oplus \Lie (H^0)\by{u} V)\in \mathrm{EHS}_{1}$ since {\it  via}\, (\ref{dia.fhs}) the following
\begin{eqnarray}\label{dia.ehs2}
\xymatrix{
H_\C\ar[r]^{\rho} & H_{\sC}\oplus \Lie (H^0)\ar[r]^{\pi}\ar[d]^u & H_{\C}\ar[d]^c \\
& V \ar[r]^{\pi_0}&  H_{\sC}/F^0_{Hodge}
}
\end{eqnarray}
commutes. The functor form  $ \mathrm{FHS}_{1}^s$ to  $\mathrm{EHS}_{1}$ just defined and the one in \ref{funEF} are clearly
mutually quasi-inverse.
\eprooff

\subsection{Hodge realization}
Deligne's Hodge realization in \cite{D} provides an equivalence
$T_{Hodge}\colon \M^{\fr}\by{\sim}{\rm MHS}_1^{\fr}$ between the  category of
Deligne $1$\nobd-motives and the category of torsion free objects in  ${\rm MHS}_1$. This equivalence has been generalized in \cite{BRS} to  an equivalence $\M \by{\sim}\MHS$ including torsion and in \cite{FHS}  to an equivalence $\Mlau \by{\sim}{\rm FHS}_1^{\fr}$ including additive  factors.
 We can now further extend both equivalences to our context (see  \ref{defaddspec} for notations).

\begin{propose}
There is an equivalence of categories
\begin{eqnarray*}
\Tint \colon  \Ma \longby{\sim}\FHS, \quad
               M\df  [u\colon F\to G]& \mapsto  \Tint (M)\df  (\Tint(F),\Lie(G))
\end{eqnarray*}
where $\Tint (M)_\et = T_{Hodge} (M_\et)$ and such that it induces an  equivalence of categories
\begin{eqnarray*}
\Tint^s \colon  \M^s \longby{\sim}\mathrm{EHS}_{1}, \quad
                 [u\colon F\to G]& \mapsto  (\Tint (M)_\et,  T_{\sC}(F_\et) \oplus \Lie (F^0)\to \Lie(G))
\end{eqnarray*}
and further restricts to the equivalence
\begin{eqnarray*}
T_{Hodge} \colon  \M \longby{\sim}\mathrm{MHS}_{1} \quad
                M & \mapsto  T_{Hodge}(M) .\end{eqnarray*}
\end{propose}
\proof The functor $\Tint$ on $\Mlau$ is constructed in \cite{FHS}.  Recall that for  $F$ free, the formal $k$\nobd-group $\Tint(F)$
is the product $F^0 \times \Tint(F_\et)$ where
 the  \'etale quotient $ \Tint(F_\et)$ is  the pull-back of $F_\et\to  G$ along $\exp{}\colon
\Lie(G)\to G$. Hence $\Tint(F_\et)$   is a free abelian group
extension of $F_\et$ by $H_1(G)$. The morphism $v\colon \Tint(F)\to
\Lie(G)$ is then taken as $\Lie(u^0)\colon \Lie(F^0)\to \Lie(G)$ on the  identity
component and the homomorphism obtained {\it via}\, the pull-back  construction
on $\Tint(F_\et)$.  The above definition of $\Tint(F)$ makes sense also  when $F$ is not free. One proceed then as done in \cite{BRS}, Prop.  1.5; in fact, one can check that $\Tint(F)$ and $\Lie(G)$ are  independent of the representative of
$ M $, \ie that a quasi-isomorphism $M_1\to M_2$ induces an isomorphism  $\Tint (M_1)\cong \Tint (M_2)$.
\eprooff

\subsection{Sharp envelope}
The sharp $\G_a$-extension of  (effective) $1$\nobd-motives, has its
counterpart in the category of formal Hodge structures. Define the
functor
\begin{eqnarray}\label{def.hsharp}
(\ \ )^\sharp: \FHS\to \FHS
\quad  (H,V) \mapsto (H,V^\sharp)
\end{eqnarray}
as follows. Let $H = H_{\Z}\times H^0\longby{v}V$. Denote $\bar v \colon H \to V/V^0$,  $c\colon  H_{\Z} \to H_{\C}$, $\iota :H^0\to \Lie (H^0)$ the identity map on Lie algebras, $\bar v_{\C} : H_{\C}\to V/V^0$ and $\bar v^0\colon H^0\to V/V^0$. Define first
\begin{itemize}
\item $(V/V^0)^{\sharp }\df H_{\C}\oplus \Lie( H^0 )$;\\
\item $\bar v^\sharp\df (c , \iota ) \colon H_{\Z}\times H^0\to (V/V^0)^{\sharp }$;\\
\item  $(V/V^0)^{\sharp 0 } \df \ker(\bar v_{\C},\bar v^0)$.
\end{itemize}
Hence we have a diagram:
\begin{eqnarray}\label{dia.vsharp}
\xymatrix {
0 \ar[r] & F^0_{Hodge}\ar[r]\ar[d]&  H_{\C} \ar[d]^{\bar v_{\C}^\sharp}\ar[r]^{\bar v_{\C}} &   V/V^0 \ar@2{-}[d] \ar[r]&  0\\
0  \ar[r]&  (V/V^0)^{\sharp 0} \ar[d] \ar[r]& (V/V^0)^{\sharp }  \ar[r]^{(\bar v_{\C},\bar v^0)} \ar[d]&  V/V^0 \ar[r] &
0 \\
 & \Lie(H^0)\ar@2{-}[r]& \Lie(H^0)\ar@/^/[u]^{\bar v^{\sharp 0}}\ar@{.>}[ur]_{\bar v^0} &&
}\end{eqnarray}
where the vertical sequence in the middle is canonically split by $\bar v^{\sharp 0}\df (\id\oplus 0)$. Define now $V^\sharp$ by pull-back as follows:
\begin{itemize}
\item $V^{\sharp }\df V\times_{V/V^0} (V/V^0)^{\sharp }$;\\
\item $v^\sharp\df (v_{\sZ}^{\sharp},v^{\sharp 0})\colon  H_{\Z} \times H^0\to V^{\sharp } $  induced by $\bar v^\sharp$ and $v$\\
\item $V^{\sharp 0}\df \ker(V^{\sharp}\to (V/V^0)^{\sharp }\to V/V^0)$.\\
\end{itemize}
Actually $V^\sharp$ fits in the following diagram of $\C$-vector spaces
\begin{eqnarray*}
\xymatrix{
0 \ar[r] & F^0_{Hodge}\ar[r]\ar[d]&   H_{\C}\times_{V/V^0} V \ar[d]^\alpha \ar[r]&  V\ar@2{-}[d] \ar[r]&  0\\
0 \ar[r]&  (V/V^0)^{\sharp 0}\ar[d] \ar[r]&  V^\sharp \ar[r]^u\ar[d]&  V  \ar[r]&
0 \\
 & \Lie(H^0)\ar@2{-}[r]& \Lie(H^0)\ar@/^/[u]^{v^{\sharp 0}} \ar@{.>}[ur]_{v^0}&&
}\end{eqnarray*}
with the canonical splitting  $v^{\sharp 0}$ of the vertical sequence in the middle. Note that the morphism $v_{\Z}^\sharp\colon H_{\Z}\to V^\sharp$ is the composition of $(c ,v_{\Z})\colon H_{\Z} \to   H_{\C}\times_{V/V^0} V$ with
$\alpha$  and the commutativity of \ref{dia.fhs} holds by construction.

\begin{remark} Note that, for  $(H, V)\in {\rm FHS}_1$, the sharp envelope $(H, V)^{\sharp}\in {\rm FHS}_1$ is such that the canonical map $v_\C^\sharp \colon H_{\sC}\to V^{\sharp}$ induced by $v_\Z^\sharp$ has a splitting $\pi:V^{\sharp}\to H_{\sC}$  induced by  $v^{\sharp 0}$.
Observe that if $(H,V)$ is \'etale, \ie $H^0=V^0=0$, we then have $(H,V)\cong (H_\Z ,H_\C/F^0_{Hodge})$ is determined by the mixed Hodge structure (\cf \cite{FHS}).  We moreover get $ (H,V)^{\sharp}\cong  (H_\Z
,H_\C/F^0_{Hodge})^\sharp = (H_\Z , H_\C)  $.
\end{remark}

\begin{lemma}\label{split}
$(H, V)\in {\rm FHS}_1^s$ if and only  if $(H, V/V^0)^{\sharp}\in {\rm FHS}_1^s$ if and only if  the splitting $\bar v^{\sharp 0}$ induces a splitting of the left most vertical sequence in \eqref{dia.vsharp}. \end{lemma}
\begin{proof} By diagram \eqref{dia.vsharp} chase, \ie $\bar v^0 : H^0\longby{\rm zero} V/V^0$ if and only if $\bar v^{\sharp 0} ( H^0)\subseteq (V/V^0)^{\sharp 0}$. 
\end{proof}

\begin{remark} Since the sharp envelope of a special formal Hodge structure is still special, 
we can consider {\it via}\, \ref{enrich} the sharp envelope on the category of enriched Hodge structrures (see  \S \ref{sec.ehs}): it  is the functor
\[(\ \ )^\sharp_e\colon {\mathrm {EHS}}_1\to {\mathrm {EHS}}_1, \quad (H_{\et}, U\to V) \mapsto (H_{\et}, U\to V^{\sharp}). \]
Here $H_{\et}$ and $U = H_{\C}\times \Lie (H^0)$ correspond to $(H, V)\in {\rm FHS}_1^s$.
\end{remark}

We further obtain a more sophisticated functor \begin{equation}\label{spsharp}
( \ \  )^{\sharp}_s : {\rm FHS}_1^s\to {\rm EHS}_1\ \ 
(H, V)\mapsto (H, V)^{\sharp}_s  \df (H_{\et}, V^{\sharp}\to V)
\end{equation}
where  $(H, V)^{\sharp}_s$ is the enriched Hodge structure associated to $\overrightarrow{(H,V)}_0$ in \eqref{overtofhs}. 
In fact, given a special structure along with its sharp envelope we are just left to get the splitting $\pi$ fitting in a commutative diagram
\[\begin{CD}
H_\C@>{v_C^\sharp}>> V^{\sharp}@>{\pi}>> H_\C \\
&& @V{u}VV @V{c}VV\\
&& V @>{\pi_0}>> H_\C/F^0_{Hodge}
\end{CD}\]
where $u: V^{\sharp}\onto V$ is the canonical projection and $\pi_0$ is the projection induced by $\sigma^{-1}: V/V^0\cong H_\C/F^0_{Hodge}$. Note that by \ref{split} and the construction of $V^{\sharp}$ we then get a natural splitting of all extensions as follows
\[
\xymatrix{
0 \ar[r] &H_{\sC}\times_{V/V^0} V\ar[r]_{\alpha}\ar[d]& \ar@/^1.3pc/[l]^{}   V^{\sharp}\ar[d]\ar[r]\ar@{.>}[dl]^{\pi}& \ar@/^1.3pc/[l]^{v^{\sharp 0}}    H^0\ar@2{-}[d] \ar[r]&  0\\
0 \ar[r]&  H_{\sC} \ar[r]& \ar@/^1.3pc/[l]^{}  (V/V^0)^{\sharp }  \ar[r]& \ar@/^1.3pc/[l]^{} H^0\ar[r]&
0\\
0 \ar[r]&  F^0_{Hodge} \ar[r]\ar[u]& \ar@/^1.3pc/[l]^{}  (V/V^0)^{\sharp 0 }  \ar[r]\ar[u]& \ar@/^1.3pc/[l]^{} H^0\ar@2{-}[u]\ar[r]&
0}\]
We obtain the claimed commutativity by diagram chase just considering that 
all these splittings are compatible. 

\begin{remark} Fix $H_{\et}=(H_{\Z},W_*,F^*_{Hodge})\in {\rm MHS}_1$ and $\pi_0 : V\to  H_\C/F^0_{Hodge}$. Then any $(H_{\et}, U\to V)\in {\rm EHS}_1$ is clearly mapped to $(H_{\et}, H_{\C}^{\sharp}\to V)$ where $H_{\C}^{\sharp}$ is just the pull-back of $H_{\C}$ along $\pi_0$. Actually, for any $(H_{\et}\times H^0, V)\in {\rm FHS}_1$ with  $V^0=\ker (\pi_0)$, we have that $(H_{\et}, H_{\C}^{\sharp})$ is the sharp envelope of $(H_{\et}, V)$ since $H_{\C}^{\sharp}\cong H_{\sC}\times_{V/V^0} V$ and $H_{\C}^{\sharp}\into V^{\sharp}$ (in particular, any $U\to V$ as above lifts to $V^{\sharp}$).
\end{remark}

Now, in \ref{sharpmot} we associated  to any effective $1$\nobd-motive $M$ the sharp extension $M^\sharp$ and  the sharp de Rham realization $T_\sharp(M)=\Lie(G^\sharp)$. Moreover, we can apply $\Tint$ to $M^\sharp$ so that we have a diagram
\[
\xymatrix{
\Ma\ar[d]_{(\ \ )^\sharp} \ar[r]^{\Tint}& \FHS\ar[d]^{(\ \ )^\sharp} \\
\Ma \ar[r]^{\Tint}& \FHS }\]

\begin{thm}\label{compthm}
Let $M=[u\colon F\to G]$ be a $\C$-$1$\nobd-motive.  The
diagram above commutes (up to isomorphisms) so that
\begin{eqnarray*}
\Tint(M)^\sharp\cong \Tint(M^\sharp)=(\Tint(F),T_\sharp(M)).
\end{eqnarray*}
\end{thm}
\proof It is sufficient to prove the commutativity on the categories
$\Mlau$  and  $\FHSf$ (\cf \ref{rem.sharp}).
Let then $M=[F\to G]$ be a Laumon $1$\nobd-motive.
We may assume  $V(G)=0=V^0$. Indeed $M^\sharp$ is obtained {\it via}\,
pull-back from $(M_\times)^\sharp=M_\times^\natural$ (\cf
\ref{lem.sharp}) and $(H,V)^\sharp$ is defined as the pull-back
of $(H,V/V^0)^\sharp$. Let then
$M=M_\times=[F\to G_\times]$ and  $M^\sharp=[F\to
G_\times^\natural]$ its universal extension. We have
\[\Tint(M)=(H,V)=(F^0\times \Tint(F_\et) , \Lie(G_\times)) \text{  with  }
V^0=0, \]
\[ \Tint(M^\sharp)=(F^0\times \Tint(F_\et),  \Lie(G_\times^\natural))\]
and  $\Tint(F_\et)\to \Lie(G_\times)$ that factors through
$\Lie(G_\times^\natural )$ because $G_\times^\natural$ is extension of
$G_\times$ by a vector group.
Denote by $[F_{\et}\to G_\et^\natural]$ the universal extension of
$M_\et$. From  \ref{rem.fbasicsh} we get a push-out diagram
\[
\xymatrix{
0 \ar[r] & \Ext(M_\et,\G_a)^\vee \ar[r]\ar[d]& G_\et^\natural  \ar[d]  \ar[r]& G_\times \ar@2{-}[d] \ar[r]&  0\\
0 \ar[r]&  \Ext(M_\times,\G_a)^\vee \ar[d] \ar[r]&  G_\times^\natural  \ar[r]\ar[d]&  G_\times \ar[r]&
0 \\
 & \Lie( {F}^0 )\ar@2{-}[r]& \Lie( {F}^0 )&&
}
\]
Consider the associated diagram of Lie algebras and compare it with (\ref{dia.vsharp}).
Recalling that $\Lie(G_\et^\natural)\cong H_\C$, $V=\Lie(G_\times)$,
$V^0=0$, $\Lie(H^0)=\Lie(F^0)$ we deduce that $V^\sharp\cong\Lie(G_\times^\natural)=T_\sharp(M_\times)$.
\eprooff

\begin{remark}
The previous theorem generalizes to (effective) Laumon  $\C$\nobd-$1$\nobd-motives
the fact that $T_{\dR}(M)\cong T_{\Z}(M)\otimes_\Z \C$ for Deligne  $1$\nobd-motives.
\end{remark}

\section{Duality on sharp de Rham realizations}

Let $M$ be  a Laumon
$1$\nobd-motive,  $M^*$ its dual and  ${\mathcal
P}\in \Biext (M, M^*;\G_m)$ the Poincar\'e biextension. We look for a canonical duality between $T_\sharp(M)$ and $T_\sharp(M^*)$ that generalizes Deligne's
construction in \cite{D} 10.2.7.3. In order to do this we need to
introduce a ``canonical'' connection on the biextension   ${\mathcal
P}^\sharp$ of $M^\sharp $ and $M^{*\sharp}$ by $\G_m$ given by the 
pull-back of ${\mathcal P}$.

\subsection{$\natural$\nobd-structures} Let $N=[F\to E]$ be an extension of $M$ by an
algebraic $k$\nobd-group $H$.
A \emph{$\natural$\nobd-structure} on $N$ is a
$\natural$\nobd-structure on $E$ as in \cite{D}, 10.2.7.1.
A \emph{strong $\natural$\nobd-structure} on
$N$ is a $\natural$\nobd-structure on $E$ such that its pull-back to
$F$ is trivial; if $F=F_\et$ any $\natural$\nobd-structure is strong.
We can characterize $G^\sharp$ as follows:

\begin{propose}\label{pro.reprGSh}
Let $M$ be a Laumon $1$\nobd-motive. Then
$G^\sharp$ is the group scheme that represents the pre-sheaf for the flat site
on $k$
\begin{equation*}
 {\mathcal F} \colon S/k \leadsto  \left\{
 \begin{array}{ll}
 (g,\nabla),& g\in G(S), \nabla\text{ a  }
\natural\text{-structure on the extension } {\mathcal P}_g \text{ of }\\
&M' \text{ by } \G_{m, S} \text{ induced by } g
\end{array}
\right\}.
\end{equation*}
\end{propose}
\proof
By (\ref{lem.Gsh}) we know that $G^\sharp$ is extension of ${\bm \omega}_{L'}$   by $A^\natural\times_A G$  and one proceeds as in \cite{BER} 3.10.
\eprooff

\subsection{The canonical connection}\label{canconn} The identity
on $G^\sharp$ provides {\it via}\, the functor $\mathcal F$ a pair  $(\rho,
\nabla_2^\sharp)$ where $\rho\colon G^\sharp\to G$ is the usual
projection and $\nabla_2^\sharp$ is a $\natural$-structure on
${\mathcal P}_\rho$ the pull-back of ${\mathcal P}$ to
$G^\sharp\times G'$ viewed as $\G_m$-extension of $G^\sharp$ over
$G'$. The same the identity on $G^{\prime\sharp}$ provides a pair
$(\rho',\nabla_1^\sharp)$ where $\nabla_1^\sharp$ is a a
$\natural$-structure on ${\mathcal P}_{\rho'}$. As ${\mathcal
P}^\sharp$ is the pull-back of ${\mathcal P}_\rho$ {\it via}\, $\rho'$  as
well as the pull-back of ${\mathcal P}_{\rho'}$ {\it via}\, $\rho$  we
define the \emph{canonical connection} $\nabla^\sharp$ on ${\mathcal
P}^\sharp$ as the sum of the (pull-back) of the connections $\nabla_i^\sharp$.
If $M$ is a Deligne $1$\nobd-motive, $\nabla^\sharp$ is the unique
$\natural$-structure on ${\mathcal P}^\sharp=P^\natural$ in \cite{D} 10.2.7.4.

\begin{example}
Let $F^0=\Spf(k[[x]])$, $M=F^0[1]$,  $M^*=F^{0*}=\Spec(k[y])$
 and $\mathcal P$ the Poincar\'e biextension. It is the trivial
$\G_m$-torsor on $F^{0*}$ together with the trivialization
$\sigma\colon F^0\otimes F^{0*}\to \G_m $ induced by Cartier duality. The
pull-back  $\mathcal P^\sharp$ of $\mathcal P$ to $([F^0 \to {\bm
\omega}_{F^{0*}}], F^{0*})$ is the trivial biextension on $({\bm
\omega}_{F^{0*}},  F^{0*} )$ together with the trivialization $\sigma$. 
The connection $\nabla_2^\sharp$ of the
trivial $\G_m$\nobd-extension of $F^{0*}$ over ${\bm
\omega}_{F^{0*}}$ is given by the invariant differential of $F^{0*}$
over ${\bm \omega}_{F^{0*}}$ associated to the identity map on ${\bm
\omega}_{F^{0*}}$; hence $xdy$. The connection $\nabla_1^\sharp$ is
associated to an invariant differentials of the $0$ group  over
$F^{0*}$ hence is trivial. In particular $\nabla^\sharp$ is
associated to $xdy$ on $  {\bm \omega}_{F^{0*}}\times_k F^{0*}$.

Observe that also $xdy+ydx$ provides a bi-invariant connection
on ${\mathcal P}^\sharp$ different from the canonical one.
Hence we can not expect a uniqueness result as in
\cite{D}, 10.2.7.4, for the (weak)   $\natural$-structures.
\end{example}

\subsection{Deligne's pairing} 
Consider the canonical connection $\nabla^\sharp$ on
${\mathcal P}^\sharp$ defined in \ref{canconn}. Its curvature is an
invariant $2$\nobd-form on $G^\sharp\times G^{\prime\sharp}$; hence
it gives an alternating pairing $R$ on
\begin{eqnarray*}
\Lie(G^\sharp\times G^{\prime\sharp})=
\Lie({G^\sharp})\oplus \Lie({G^{\prime\sharp}})= T_\sharp(M)\oplus
T_\sharp(M^*)
\end{eqnarray*}
 with values in $k$. As the restrictions of $R$ to
$\Lie(G^\sharp)$  and  $\Lie({G^{\prime\sharp}})$  are trivial it
holds
\begin{eqnarray*}
R(g_1+g_2,g_1'+g_2')=\Phi(g_1,g_2')-\Phi(g_2,g_1')
\end{eqnarray*}
with
\begin{eqnarray}\label{pairing}
\Phi\colon T_\sharp(M)\otimes T_\sharp(M^*)\to k.
\end{eqnarray}

If $M$ is a Deligne $1$\nobd-motive, the pairing
above coincides with the one in \cite{D} 10.2.7.

We will see that $\Phi$ is perfect following the proof
in \cite{BER}, \S 4, for the classical case of Deligne $1$\nobd-motives.

Recall  the  extensions  in (\ref{diasharp}) for $M$ and $M^*$:
\begin{eqnarray}\label{seq.wEMG}
\xymatrix@C=17pt{
 &0\ar[r]&{\bm \omega}_{G'}\ar[r]^i & M^\sharp \ar[r]^\rho& M\ar[r]&
 0,\quad
 0\ar[r]&{\bm \omega}_{G}\ar[r]^{i'} & M^{*\sharp} \ar[r]^{\rho'}&  M^*\ar[r]&
 0,
}
\end{eqnarray}
We denoted by  $({\mathcal P}_{\rho},\nabla_2^\sharp)$  the
 $\natural$\nobd-extension of $M'$ by the
multiplicative group over $G^\sharp$ that corresponds to the
identity map on $G^\sharp$ {\it via}\, the functor $\mathcal F$ in
\ref{pro.reprGSh}. Similarly for  $({\mathcal
P}_{\rho'},\nabla_1^\sharp)$.

\begin{lemma}\label{lem.alpha} Let
$\alpha_{G'}$ be the invariant differential of $G'$ over ${\bm
\omega}_{G'}$ that corresponds to the identity map on ${\bm
\omega}_{G'}$. The restriction of $({\mathcal P}_{\rho
},\nabla_2^\sharp)$ to ${\bm \omega}_{G'}$ {\it via}\, $i\colon {\bm
\omega}_{G'}\to G^\sharp$ in (\ref{seq.wEMG}) is isomorphic to the
trivial extension of $ M^*$ by the multiplicative group over ${\bm
\omega}_{G'}$ equipped with the connection associated to
$\alpha_{G'}$.
\end{lemma}
\proof See \cite{BER}, 4.1. \eprooff

Changing the role of $M$ and $M^*$, denote by $\alpha_{G}$
 the invariant differential of $G$ over ${\bm \omega}_{G}$
that corresponds to the identity map on ${\bm \omega}_{G}$. The
restriction  of $({\mathcal P}_{\rho'},\nabla_1^\sharp)$ to ${\bm
\omega}_{G}$ is isomorphic to the trivial extension of $G'$ by the
multiplicative group over ${\bm \omega}_{G}$ equipped with the
connection associated to $\alpha_{G}$. From \cite{BER}, 4.2, we know
that
\begin{lemma}\label{lem.dalpha}
The curvature of $\alpha_G$ provides a perfect pairing
\begin{eqnarray*}
d\alpha_G\colon~{ \omega}_{G}\otimes \Lie({G}) \longrightarrow k
\end{eqnarray*}
that is the usual duality.
\end{lemma}
Hence the proof of Theorem 4.3 in \emph{loc. cit.} works the same
and we get
\begin{thm}\label{thm.perfM}
Let $M$ be a free $k$-$1$\nobd-motive. The pairing $\Phi$ in
$(\ref{pairing})$ is perfect. Moreover it fits in a diagram
\begin{eqnarray*}
\xymatrix @C-=1pt{
{\omega}_{G'}\ar[d]^{\iota}&\otimes& \Lie({G'})\ar[rr]&\quad & k\\
\Phi\colon \Lie({G^\natural})\ar[d]^{g}& \otimes &
 \Lie(G^{\prime\sharp})\ar[u]^{g'}\ar[rr]& &k\\
\Lie({G})&\otimes&{ \omega}_{G}\ar[rr]\ar[u]^{\iota'}& &k}
\end{eqnarray*}
where the vertical homomorphisms come from
(\ref{seq.wEMG}) and the upper (resp. lower) pairing is the usual duality between the Lie algebra of $G'$ (resp. $G$) and the $k$\nobd-vector space of  invariant
differentials of $G'$ (resp. $G$).
\end{thm}

\section{Sharp de Rham cohomology}
We describe $H^1_{\sharp-\dR}(X)\df T_{\sharp}(\Pic^+_a(X))$ in some meaningful cases, \ie  when $X$ is proper or is a smooth algebraic $k$-scheme. Here $\Pic_a^+ (X)$ and its Cartier dual  $\Alb_a^-(X)$ are the Laumon 1-motives of the algebraic $k$-scheme $X$  constructed in \cite{SAP}, whence $\Pic_a^+ (X)_{\et} =\Pic^+ (X)$ and $\Alb_a^-(X)_{\et} = \Alb^-(X)$ were introduced in \cite{BVS}. By construction, see~\ref{sharpmot}, we then have that  $H^1_{\sharp-\dR}(X)$ is sitting in an extension
$$0\to H^1_{\dR}(X) \to H^1_{\sharp-\dR}(X) /V(\Pic) \to V (\Alb) \to 0$$
where we  have set for $\Pic_a^+ (X)= [F \to G]$
\begin{itemize}
\item $V (\Pic)\df $ the additive part, \ie (the Lie algebra of) the vector group $V(G)$ given by the maximal additive subgroup of $\Pic^{0}(\bar X)$ for a suitable (singular) compactification $\bar X$ of $X$ 
\item $V(\Alb ) \df$ the infinitesimal part, \ie the Lie algebra $\Lie F^0$ that is just the dual of the corresponding Faltings-W\"ustholz vector group in the Albanese $\Alb_a^-(X)$.
\end{itemize}
Here $V (\Pic) =0$ if $X$ is smooth and $V (\Alb) =0$ if $X$ is proper over $k$.
Moreover, for $X$ smooth we have that $F_{\et} =\Div_{Y}^0(\bar X)$ where $\bar X$ is a smooth proper compactification, $Y =\bar X -X$ is a normal crossing divisor, and $\Lie F^0$ is the $k$-vector space $\ker( H^1(\bar X, \cO_{\bar X})\to H^1(X,\cO_X))$, \ie is $\Gamma (X,\Omega^1_X)_{\rm closed}/d(\Gamma (X, \cO_X))$ divided out by $\Gamma (\bar X, \Omega^1_{\bar X}(\log Y))$. For $X$ proper over $k$ it holds $F =0$ and $G = \Pic^0 (X)$ is the connected algebraic group given by the identity component of the representable fppf-sheaf $\Pic_{X/k}$ (see \cite{SAP} for more details).

\subsection{Sharp extension of $\Pic^+_a(X)$}\label{picsharp} We compute the sharp $\G_a$-extension of $\Pic_a^+ (X)$ for $X$ proper or smooth. 

For $X$ proper and $\Xs\to  X$ a smooth proper hypercovering  we obtain the semi-abelin quotient $\Pic^{0}(X)/ V (\Pic) =  \bPic^{0}( \Xs)$ by \cite{BVS}, Lemma 5.1.2. In \loccit we also introduced the algebraic group $\bPic^{\natural}(X_{\d})$ given by 
isomorphism classes of triples $(\cL,\nabla, \alpha)$ consisting of
an invertible sheaf $\cL$ on $X_0$, with an integrable connection
$\nabla$, and an isomorphism $\alpha :d_0^*(\cL,\nabla)\by{\cong}
d_1^*(\cL,\nabla)$ satisfying the cocycle condition (here $d_0, d_1:X_1\to X_0$ are the face maps). There is a functorial isomorphism
\begin{equation}\label{picnat}\bPic^{\natural}(X_{\d})\cong
\HH^1(X_{\d},\cO^*_{X_{\d}}\by{\rm dlog}\Omega^1_{X_{\d}})
\end{equation}
Define the group scheme $\Pic^{\sharp} (X)\df
\bPic^{\natural} (\Xs) \times_{\bPic (\Xs)}\Pic (X)$ by pull-back.

\begin{lemma}
If $X$ is proper then $\Pic_a^+ (X)^{\sharp} = [0\to \Pic^{\sharp, 0} (X) ]$.
\end{lemma}
\proof
Since $X_{\d}$ is smooth and proper over $k$, the
semi-abelian variety $\bPic^0(X_{\d})$ is mapped to zero in
$\HH^1(X_{\d},\Omega^1_{X_{\d}})$ {\it via}\, \eqref{picnat}.
 We then have that $\bPic^{\natural, 0}(X_{\d})$ is an extension of $\bPic^0(X_{\d})$ by  $\HH^0(X_{\d},\Omega^1_{X_{\d}})$, \ie is the pull back of the inclusion $\bPic^0\into \bPic$. This extension is the universal $\G_a$-extension of the semi-abelian scheme $\bPic^0(X_{\d})$ by \cite{BVS}, Lemma 4.5.2. Since $V(\Pic) = \ker (\Pic^{\sharp, 0} (X) \onto\bPic^{\natural, 0} (\Xs))$ we then get the following pullback diagram which, by \eqref{diasharp}, proves the assertion
\begin{equation} \label{picsh}
\xymatrix{
0 \ar[r] & \HH^0(\Xs, \Omega^1_{\Xs})\ar[r] \ar@2{-}[d]  & \Pic^{\sharp, 0} (X)  \ar[d]  \ar[r]& \Pic^{0} (X)  \ar[d] \ar[r]&  0\\
0 \ar[r] & \HH^0(\Xs, \Omega^1_{\Xs}) \ar[r]&  \bPic^{\natural, 0} (\Xs) \ar[r]&  \bPic^{0} (\Xs)  \ar[r]& 0 }
\end{equation}
\eprooff

For $X$ smooth recall the algebraic group $\Pic^{\natural-\log}(X)$ given by isomorphism classes of pairs $(\cL ,\nabla^{\log})$ where $\cL$ is a line bundle on $\bar X$ and $\nabla^{\log}$ is an integrable connection on $\cL$ with log poles along $Y$. In \cite{BVS}, Lemma 2.6.2, we have seen that $\Pic^+_a(X)_{\et}^{\natural} = [\Div_Y^0(\bar X)\to \Pic^{\natural-\log, 0}(X)]$. Actually, the universal $\G_a$-extension of $\Pic^+_a(X)$ exists (\cf the remark \ref{rem.fbasicsh}):
\begin{lemma} 
If $X$ is smooth then
$\Pic_a^+ (X)^{\natural} = [F\to \Pic^{\natural-\log, 0}(X) + \Lie F^0]$.
\end{lemma}

\subsection{Sharp de Rham over $\C$}\label{drsharp} 
For $X$ over $\C$ we also have
\[H^1_{\sharp-\dR}(X) = T_{\oint}(\Pic^+_a(X))^{\sharp}\]
 by \ref{compthm}.  Let $X$ be a proper $\C$-scheme and $\Xs$ a smooth proper hypercovering as above. In this case, passing to the Lie algebra  $\Lie \bPic^{\natural, 0} (\Xs) = H^1_{\dR}(X) \cong   \HH^1(\Xs,\C)  \cong H^1(X,\C)$ by the (simplicial) de Rham theorem and cohomological descent for the analytic topology, \cf \cite{BVS}, remark 2.6.3. We then obtain, \cf \ref{spsharp}, an enriched Hodge structure {\it via}\, the following diagram of the Lie algebras of \eqref{picsh}
\begin{equation}\label{sharpicsq}
\xymatrix{
0 \ar[r] & \HH^0(\Xs, \Omega^1_{\Xs})\ar[r] \ar@2{-}[d]  & H^1_{\sharp-\dR} (X)  \ar[d]^\pi  \ar[r]^u & H^1(X,\cO_X) \ar[d]^{\pi_0}\ar[r]&  0\\
0 \ar[r] & \HH^0(\Xs, \Omega^1_{\Xs}) \ar[r]& H^1(X,\C)\ar[r]^c&  H^1(\Xs,\cO_{\Xs})  \ar[r]& 0 }
\end{equation}
where $H^1(X,\C)/F^1_{Hodge} =H^1(\Xs,\cO_{\Xs})  $, $\HH^0(\Xs, \Omega^1_{\Xs}) =F^1_{Hodge}$ and $V(\Pic) =\ker (H^1(X,\cO_X)\onto H^1(\Xs,\cO_{\Xs}))$. If $X$ is smooth we have $\Lie\Pic^{\natural-\log, 0}(X) = H^1(X, \C)$ by \cite{BVS}, 2.6.4. We then obtain:

\begin{propose} 
For $X$ over $\C$ we have $H^1_{\sharp-\dR} (X)\cong H^1(X,\C ) \oplus V $ where $V = V(\Pic)$ if $X$ is proper and $V=V(\Alb) $ if $X$ is smooth.
\end{propose}
\begin{proof} It follows from the previous lemmas in \ref{picsharp} and the above discussion.\end{proof}

\begin{cor} 
If $X$ is a proper $\C$-scheme and $H^1(X,\Z)=0$ then $H^1_{\sharp-\dR} (X) = H^1(X,\cO_X)$.
\end{cor}

Let $X$ be now a proper (reduced) variety over $\C$ and consider, following \cite{BSE}, the naive analytic de Rham complex $\Omega^{\d}_X$ on $X$ itself. The resulting cohomology $\HH^n(X,\Omega^{\d}_X)$ is considered in \cite{BSE} as part of one possible enriched Hodge structure associated to $X$ and $n$. Actually for $'\Omega^{\d}_X\df
\Omega^{\d}_X/F^{\dim (X) +1}\Omega^{\d}_X$ and
$''\Omega^{\d}_X := \Omega^{\d}_X/\text{tors}$ we have
$\Omega^{\d}_X \onto\ '\Omega^{\d}_X\onto\
''\Omega^{\d}_X$. Moreover $$\HH^n(X,\Omega^{\d}_X) \to \HH^n(X,\
'\Omega^{\d}_X) \to \HH^n(X,\ ''\Omega^{\d}_X)
\to \HH^n(\Xs,\Omega^{\d}_{\Xs}) \cong H^n(X,\C)
$$ yields three different enriched Hodge structures associated to $X$ and $n$,
see \cite{BSE}, 2.1 and 2.2. For $n=1$ we want to compare them with  $H^1_{\sharp-\dR} (X)$.
We clearly have, by construction, a commutative diagram with exact rows
\begin{equation}\label{comp}
\xymatrix{
0 \ar[r] & \HH^0(X,\ ''\Omega^{\d\geq 1}_{X})\ar[r] \ar[d]  &  \HH^1(X,\ ''\Omega^{\d}_X)\ar[d]  \ar[r]& H^1(X,\cO_X) \ar[d] &   \\
0 \ar[r] & \HH^0(\Xs, \Omega^1_{\Xs}) \ar[r]& H^1(X,\C)\ar[r]&  H^1(\Xs,\cO_{\Xs})  \ar[r]& 0 }
\end{equation}
providing a canonical comparison map 
\[
\rho'':\HH^1(X,\ ''\Omega^{\d}_X)\to H^1_{\sharp-\dR} (X)
\]
and similarly by composition $\rho':\HH^1(X,\ '\Omega^{\d}_X)\to H^1_{\sharp-\dR} (X)$ and $\rho:\HH^1(X,\Omega^{\d}_X)\to H^1_{\sharp-\dR} (X)$. In fact, just observe that $H^1_{\sharp-\dR} (X)$ is given by the pull-back diagram \eqref{sharpicsq} and therefore {\it any}\, such enriched Hodge structure associated to $X$ and $n=1$ maps to it.

\begin{propose} 
Let $\rho^{\dag}$ denote any comparison map $\rho, \rho'$  or $\rho''$ corresponding to ${}^{\dag}\Omega^{\d}_{X}$ which denotes  $\ \Omega^{\d}_{X}, \ '\Omega^{\d}_{X}$ and $\ ''\Omega^{\d}_{X}$ respectively. The decorated comparison map $\rho^{\dag}$ is 
\begin{enumerate}
\item[(i)] surjective if and only if the boundary map 
$H^1(X,\cO_X)\to \HH^1(X,{}^{\dag}\Omega^{\d\geq 1}_{X})$ is zero and $\HH^0(X,{}^{\dag}\Omega^{\d\geq 1}_{X})\onto \HH^0(\Xs, \Omega^1_{\Xs}) $ is surjective;
\item[(ii)] injective if and only if the map $\HH^0(X,{}^{\dag}\Omega^{\d\geq 1}_{X})\into\HH^0(\Xs, \Omega^1_{\Xs}) $ is an inclusion.
\end{enumerate}
The map $\rho^{\dag}$ is then an isomorphism if and only if both conditions hold.
\end{propose}
\proof
Comparing the diagram \eqref{sharpicsq} and the decorated version of \eqref{comp} observe that the top exact sequence in the latter continues on the right with the mentioned boundary map.
\eprooff

\begin{remark} If $\pi_1 (X) =0$ then $\rho^{\dag}:\HH^1(X,{}^{\dag}\Omega^{\d}_{X})\to H^1(X,\cO_X)=H^1_{\sharp-\dR}(X)$. 
For example, if we take the curve considered in 2.3 of \cite{BSE} then $\rho''$ is an isomorphism. Note that it seems puzzling to study the geometric meaning of these conditions. In general, we just have that $$\HH^0(\Xs, \Omega^1_{\Xs})=\ker H^0(X_0, \Omega^1_{X_0}) \longby{d_0^*- d_1^*} H^0(X_1, \Omega^1_{X_1})$$ for the components $X_0, X_1$ of $\Xs$ while $$\HH^0(X,{}^{\dag}\Omega^{\d\geq 1}_{X})= \ker H^0(X,{}^{\dag}\Omega^1_X)\to H^0(X,{}^{\dag}\Omega^2_X)$$ 
For $\dim (X)=1$,  by choosing $\Xs$ in such a way that $X_0$ is the normalization of $X$ and $X_1$ is $0$-dimensional, we have  $\HH^0(\Xs, \Omega^1_{\Xs}) = H^0(X_0, \Omega^1_{X_0}) $, $\HH^0(X,\ '\Omega^{\d\geq 1}_{X}) = H^0(X,\Omega^{1}_{X}) $ and further $\HH^0(X,\ ''\Omega^{\d\geq 1}_{X}) = H^0(X,\Omega^{1}_{X}/\text{tors})$, \eg the injectivity of $\rho'$ means that $H^0(X,\Omega^{1}_{X})$ injects into $H^0(X_0, \Omega^1_{X_0})$.
\end{remark}

\appendix
 \section{} In this section we recall some facts and results on (formal) $k$\nobd-groups needed in the paper. The characteristic of the field $k$ is zero.

 \subsection{Vector groups}
Let $\mathcal E$ be a  free $k$\nobd-module and
$\mathcal E^\vee=\Hom(\mathcal E,k)$ its dual.
Denote by $E=\Spec(k[{\mathcal E}^\vee])$ the $k$\nobd-vector group
associated to $\mathcal E$ where
\[k[{\mathcal E}^\vee]={\mathrm
{Sym}({\mathcal E}^\vee)}=k\oplus  {\mathcal E}^\vee \oplus  S^2({\mathcal
E}^\vee)\oplus \dots\]
Its completion at the origin is
$\widehat E=\Spf(k[[{\mathcal E}^\vee]])$ where $k[[{\mathcal
E}^\vee]]$ means the infinite product \[ k\times {\mathcal
E}^\vee\times  S^2({\mathcal E}^\vee)\times S^3({\mathcal
E}^\vee)\times \dots \] with the multiplication induced by that of
$\mathrm {Sym}({\mathcal E}^\vee)$.

\begin{remark}\label{rem.formlie}
Starting with $\hat E$ we can recover
$E$ {\it via}\, ${\mathcal E}=\Lie(\hat E)$ (\cf \cite{SGA3}, VII${}_B$  3.3).
\end{remark}

\subsection{On Cartier duals}\label{sec.cartierdual}
  Let $F=\Spf(A)$ be a connected formal $k$\nobd-group.
 Its Cartier dual\footnote{It is denoted by $\mathbb D(F)$ in  \cite{FO}.}
 $F^*$ is defined as $\Spec(A^*)$ with
 $A^*:=\mathrm{Hom}_{\mathrm{cont}}(A,k)$
 where $k$ is endowed with the discrete topology. For example, if
 $A=k[[x]]$, any continuous $k$\nobd-linear map $f\colon A\to k$ factors
 through $k[[x]]/(x^n)$ because $f^{-1}(0)$ has to be open in $A$.
 Set $f(1)=a_0, f(x^i)=a_i$; then $f$ is uniquely determined by the
 polynomial $\sum_i a_i(x^*)^i$. Hence $A^*=k[x^*]$. Observe that
 $x^*$ corresponds to the $k$\nobd-linear map sending $1\mapsto 0$,
 $x\mapsto 1$ and $x^i\mapsto 0$ for $i>1$. Similarly
 $k[[x_1,\dots,x_n]]^*=k[x^*_i,\dots,x_n^*]$.

The duality between $F^0$ and $F^{0*}$ provides also a duality on
 Lie algebras.

\begin{lemma}\label{lemliedual}
Let $F=\Spf(k[[x_1,\cdots,x_n]]) $ be a connected formal $k$\nobd-group
and $F^*=$ $\Spec( k[[x_1,\cdots,x_n]]^*)$ its Cartier dual.
There is a canonical duality between $\Lie(F)$ and
$\Lie(F^*)$.
\end{lemma}
\proof The  Lie algebra of $F$ corresponds to the
$k$\nobd-linear maps
\[k[[x_1,\cdots,x_n]]\to k[\epsilon]/(\epsilon^2)
\]
 such that $a\mapsto a$ for $a\in k$, $x_i\mapsto b_i\epsilon$ with
 $b_i\in k$ and $x_ix_j\mapsto 0$, hence to the $k$\nobd-linear  polynomial
 $\sum_{i=1}^n b_ix_i^*$. Recall now that
\[
F^*=\Spec(k[[x_1,\cdots,x_n]]^*)=\Spec(k[x_1^*,\dots,x_n^*]).
\]
 The Lie algebra of $F^*$ is the $k$\nobd-module of $k$\nobd-linear  polynomials
 in the $n$-variables $x_i^{**}$. Hence there is a canonical pairing
 $\Lie(F)\otimes \Lie(F^*)\to k$ sending $ x_i^*\otimes x_j^{**}$ to
 $\delta_{ij}$ that does not depend on the choice of the basis $x_i$.
 \eprooff

\subsection{Formal completion at the origin}\label{subformcompl}
Let $G$ be a connected algebraic $k$\nobd-group. The connected formal
$k$\nobd-group associated to $\Lie(G)$ is canonically isomorphic to the
formal completion at the origin of $G$. Indeed, let $x_1,\cdots,x_g$
be free generators of $\Lie(G)^\vee$ over $k$. The associated formal
$k$\nobd-group is $\Spf(k[[x_1,\cdots,x_g]])$; moreover, as
$\Lie(G)^\vee$ is canonically isomorphic to the $k$\nobd-module of
invariant differentials on $G$ and hence to $m/ m^2$ with $m$ the
maximal ideal of ${\mathcal O}_{G,e}$, we could think $\{x_i\}_i$ as
a basis of the $k$\nobd-module $ m/ m^2$. Now, the formal completion at
the origin of $G$ is the formal spectrum of
\[\liminv n {\mathcal O}_{G,e}/m^n=k\times  m/
m^2\times m^2/ m^3\times \cdots=k[[x_1,\dots,x_g]].
\]

As a consequence of \ref{lemliedual} we get then

\begin{lemma}\label{lem.conno}
Let $G$ be an  algebraic $k$\nobd-group.
The Cartier dual of its formal completion $\hat G$ is canonically
isomorphic to  ${\bm \omega}_{G}=\Spec(k[\Lie(G)])$.
\end{lemma}

\subsection{Homomorphisms and extensions}
We defined ``strongly exact'' sequences in $\eMa$ as exact sequences 
of complexes in ${\bf Ab}/k$. It is immediate to prove the following:
\begin{lemma}\label{lem.ext}
Il $M$ is an effective $k$-$1$-motive and  $E$ is any sheaf in 
${\bf Ab}/k$ the morphism
\[ {\mathrm{Ext}}_{C^{[-1,0]}({\bf Ab}/k)}(M,E)\longrightarrow 
{\mathrm {Hom}}_{D^b({\bf Ab}/k)  }(M,E[1])\]
is an isomorphism.  
\end{lemma}

Now, $\eMa$ is an exact subcategory of $\Ma$ and
strongly exact sequences of effective $1$-motives are exact in $\Ma$. The converse is not true in general. However, any exact sequence in $\Ma$ can be represented by a strongly exact sequence (\ref{cnsdual}, \cf \cite{BRS} for the classical case). Furthermore,  $\G_a$-extensions  of $1$-motives are isomorphic to strongly exact extensions.

\begin{propose}\label{pro.ext}
Let $M$ be an effective $k$-$1$-motive and $W $ a $k$\nobd-vector group.
Any isomorphism class of extensions of $M$ by  $W$ in $\Ma$ contains a
strongly exact extension of $M$ by $W$ and the canonical map
\[ {\mathrm{Ext}}_{C^{[-1,0]}({\bf Ab}/k)}(M,W)\longrightarrow {\mathrm {Ext}}_{\Ma}(M,W) \] 
is an isomorphism.
\end{propose}
\proof
The injectivity follows immediately from the fact that any q.i. between $1$-motives 
$[F\to G_i]$, $i=1,2$, is an isomorphism.
For the surjectivity, 
let $\tilde G\to [F_1\to G_1]\by{f} M$ be effective morphisms that provide an extension in $\Ma$. 
It means that $\iota\colon \tilde G\to G_1$ is a monomorphism and 
$f$ induces  epimorphisms $f_F\colon F_1\to F$ on the formal groups 
and $f_G\colon G_1\to G$ on the algebraic
$k$-groups. Moreover, $\ker(f_F)=\ker(f_G)/\im(\iota)$.
If now $\tilde G$ is a vector group $W$, one deduces easily that $W$ is the kernel of 
the restriction of $f$ to $V(G_1)$ and that any extension $[F_1\to G_1]$ of $M$ by $W$ in $\Ma$ is isomorphic to 
the extension $[F\to G_1/\ker(f_F)]$ of $M$ by $W$ in $\eMa$. 
\eprooff

\begin{lemma}
\label{lem.homform} Let $F$ be a formal $k$-group. Then
$\Hom(F,\G_a)$ is a free $k$\nobd-module of finite rank.
\end{lemma}
 \proof
 For the connected part,  
it is sufficient to consider the case $F^0=\fagr$. Now,
 $\Hom(\fagr,\agr)=\Hom(\fagr,\fagr)=k$. For $F$ \'etale and free, $
\Hom(F,\agr)=\Lie(T')$ where $T'$ is the Cartier dual of $F$.
For $F_{\tor}$  one has $\Hom(F_{\tor},\G_a)=0$. \eprooff

\begin{lemma}\label{homtrivfg}
Let $G$ be a connected algebraic $k$\nobd-group and $F$ a formal
$k$\nobd-group. Then $\Hom(G,F)=0$.
\end{lemma}
\proof As $G$ is connected, $\Hom(G,F_\et)=0$. It remains to prove
that $\Hom(G,\fagr)=0$.  Any morphism $f\colon L\to \fagr$, with $L$ a linear $k$-group, is
trivial because $L$ is reduced.  Suppose then $G=A$ an abelian
variety and let $f\colon A\to \fagr$ be a morphism. The induced morphism
$A\to \G_a$ is trivial. Moreover for any $k$\nobd-algebra $C$,  $\fagr(C)={\mathrm{Nil}}(C)$ 
injects into $\G_a(C)=C$. Hence also $f$  is trivial.  \eprooff

\begin{lemma}\label{exttrivfg}
Let $F$ be a  formal $k$\nobd-group without torsion and $G$ an algebraic
connected $k$\nobd-group. Then
$\Ext_{{\bf Ab}/k}(G,F)=0$.
\end{lemma}
\proof
For $F=F_\et$ see \cite{RA}, 2.3.2. For $F=F^0$ we reduce to the
case $F=\fagr$. Observe that
\begin{equation*}
\Ext_{{\bf Ab}/k}(G,\fagr)=\Ext_{{\bf Ab}/k}(G,\ihom(\G_a,\G_m))=\biext{G}{\G_a}{\G_m}
\end{equation*}
where we use Cartier duality for the first isomorphism and the exact
sequence
\begin{equation*}
\Ext_{{\bf Ab}/k}(P,\ihom(Q,\G_m))\to \biext{P}{Q}{\G_m}\to
\Hom(P,\uext{Q}{\G_m})
\end{equation*}
(\cf \cite{SGA7}, VIII, 1.1.4) with $P=G$, $Q=\G_a$ for the second.
Moreover
\[\biext{G}{\G_a}{\G_m}=\biext{\G_a}{G}{\G_m}=0
\]
because of \cite{SGA7}, VII, 3.6.5 \& VIII, 4.6.
\eprooff

\begin{lemma}\label{exttrivfa}
Let $F$ be a connected formal $k$\nobd-group. Then 
$\Ext_{ {\bf Ab}/C }({F},{\G_a})=0$ for any  $k$-algebra $C$
and ${\underline{\Ext}}(F,\G_m)=0$.
\end{lemma}
\proof We  reduce to the case $F=\fagr$.  
Suppose given an extension
\[0\to \G_a\to H\to \fagr\to 0\]
over $k$. We show that it is trivial.
 Denote by $H_n$ the pull-back of $H$ to the $n$\nobd-infinitesimal 
neighborhood  $\G_{a,n}=\Spec(k[x]/(x^{n+1}))$ of
$\fagr$.  
The scheme $H_n$ is an $\G_a$-torsor over $\G_{a,n}$ and hence trivial.
In particular $H_n$ is smooth over $\G_{a,n}$. Recalling now that  $\G_{a,n-1}\to \G_{a,  n}$ is a closed immersion with square zero ideal, the lifting property of smooth morphisms permits to construct  a tower of compatible sections 
$s_n\colon \G_{a,n}\to  H_n\to H$  and hence compatible ``factor sets''
$\gamma_n \colon \G_{a,n}^2\to \G_a$ defined as  
\[\gamma_n(a,b)\df s_{2n}(a+b)-s_n(a)-s_n(b).
\]
 We may suppose that $\gamma$ is normalized, \ie $s_n(0)=0\in \G_a(k)$.
We may summarize this fact saying that we have a morphism $\gamma\colon \fagr \to \G_a$ as contravariant functor from ${\bf Aff}/k$ to the category of sets satisfying the usual properties of a  factor set. Let $P=\sum a_{ij}x_1^ix_2^j$ be the associated power series  in $k[[x_1,x_2]]$.  As $\gamma_n(0,b)=\gamma_n(a,0)=s_n(0)=0$ the polynomial  
$\gamma_n^*(x)$ (that is $P$ truncated at the $n$th powers) is divisible by $x_1x_2$ and  
$\gamma_n$ factors through $\G_{a,n}$. It provides then a ``factor set''
 $\hat\gamma\colon \fagr^2\to \fagr$ and $H$ is the push-out along $\fagr\by{\iota}\G_a$ of an extension $E$ of $\fagr$ by itself. As any extension of connected formal $k$-groups  is trivial, $E$ is trivial and hence 
 the same is $H$.

Let now $C$ be a $k$-algebra and $H$ an extension of $\fagr$ by $\G_a$ over $C$.
We can repeat the above construction getting  ``factor sets'' $\gamma_n$ over $C$
determined by a power series $P=\sum a_{ij}x_1^ix_2^j$ in $C[[t]]$. 
In order to see that $H$ is trivial, one is reduced to see that $P$ can be written as $h(x_1+x_2)-h(x_1)-h(x_2)$ for a suitable power series $h(t)=\sum b_nt^n$ with coefficients in $C$. This is possible if 
\begin{eqnarray}\label{eq.cond}\frac{a_{ij}}{{i+j\choose i}}=\frac{a_{mn}}{{m+n\choose m}}
\quad {\rm whenever}\quad i+j=m+n.\end{eqnarray}
This fact can be deduced comparing 
\[\sum a_{ij}x_1^ix_2^j \quad {\rm and} \quad \sum b_n(x_1+x_2)^n- 
\sum b_nx_1^n-\sum b_nx_2^n.\]
Using now the property of factor sets
\[\gamma(a,b)-\gamma(a, b+c)=\gamma(b,c)-\gamma(a+b,c)\]
we get 
\begin{eqnarray}\label{eq.series}
\sum a_{ij} x^iy^j       -\sum a_{ij} x^i (y+z)^j = 
\sum a_{ij} y^iz^j       -\sum a_{ij}(x+y)^i z^j,    \quad {\rm or} \nonumber \\
\sum  a_{ij} x^i[y^j-(y+z)^j]=\sum a_{ij}[y^i-(x+y)^i] z^j.
\end{eqnarray}
Now $a_{ij}{j\choose s}$ is the coefficient of the term $-x^iy^{j-s}z^s$ while
$a_{mn}{m\choose q}$ is the coefficient of the term $-x^{m-q}y^{q}z^n$.
Observe that a monomial $x^ay^bz^c$  occurs once on the left for $i=a$ and $j=b+c$ and on the right for $m=a+b$ and $n=c$. 
Considering the same monomial on the right and on the left, the indices
satisfy the following relations: $i+j=m+n$, $s=n$ and $q=j-n$.
The \ref{eq.series} implies that 
\[a_{ij}{j\choose n}=a_{mn}{m\choose j-n}.\]
As
\[{j\choose n}{i+j\choose i}={m\choose j-n}{m+n\choose m}\] 
condition \ref{eq.cond} is satisfies.

Let now $H$ be an extension of $\fagr$ by $\G_m=\Spec(k[z,z^{-1}])$ over a $k$-algebra $C$. 
Let $H_n$ be the pull-back of $H$ to $\G_{a,n}$.
\'Etale locally on $C$, $H_n$ is a trivial $\G_m$-torsor. We  may suppose that it is trivial on $C$. Proceed then as above constructing sections $s_n\colon \G_{a,n}\to H$ and normalized
 ``factor sets'' $\gamma_n$ such that $\gamma_n^*(z-1)$ is divisible by
$x_1x_2$. Hence $\gamma_n$ factors through $\Spec(C[z,z^{-1}]/(z-1)^n)$ and
we get a factor set $\gamma\colon \fagr^2\to \hat \G_m$ over $C$.
As $\hat\G_m$ is isomorphic to $\fagr$, we proceed as done for $\G_a$,  
showing that $\gamma$ is equivalent to the trivial factor set.
\eprooff

\hrulefill

 \end{document}